\documentclass[aap]{imsart}
\usepackage{amsmath,amstext}  
\usepackage{amsfonts,amssymb} 
\usepackage[latin1]{inputenc}
\usepackage[english]{babel}
\usepackage{graphics}
\usepackage{epsfig}
\usepackage{epsf}
\usepackage{here}
\usepackage{color}
\usepackage{pifont}

\startlocaldefs
\numberwithin{equation}{section}

\newcommand{\dT}{\mathbb{T}}
\newcommand{\dG}{\mathbb{G}}
\newcommand{\dF}{\mathbb{F}}
\newcommand{\dR}{\mathbb{R}}
\newcommand{\dE}{\mathbb{E}}

\newcommand{\dX}{\mathbb{X}}
\newcommand{\dU}{\mathbb{U}}
\newcommand{\cA}{\mathcal{A}}
\newcommand{\cN}{\mathcal{N}}

\newcommand{\cB}{\mathcal{B}}

\newcommand{\cO}{\mathcal{O}}
\newcommand{\cR}{\mathcal{R}}
\newcommand{\cT}{\mathcal{T}}

\newcommand{\cF}{\mathcal{F}}

\newcommand{\cV}{\mathcal{V}}
\newcommand{\cW}{\mathcal{W}}
\newcommand{\rI}{\mathrm{I}}
\newcommand{\rJ}{\mathrm{J}}
\newcommand{\cM}{<\!\!M\!\!>}
\newcommand{\veps}{\varepsilon}
\newcommand{\wh}{\widehat}

\newcommand{\pt}{\!\!\!\textbf .}
\newcommand{\ind}{\mbox{1}\kern-.25em \mbox{I}}
\font\calcal=cmsy10 scaled\magstep1
\def\build#1_#2^#3{\mathrel{\mathop{\kern 0pt#1}\limits_{#2}^{#3}}}

\def\videbox{\mathbin{\vbox{\hrule\hbox{\vrule height1ex \kern.5em
\vrule height1ex}\hrule}}}
\endlocaldefs

\begin{document}
\newtheorem{Remark}{Remark}[section]
\newtheorem{Proposition}[Remark]{Proposition}
\newtheorem{Theorem}{Theorem}[section]
\newtheorem{Lemma}[Remark]{Lemma}
\newtheorem{Corollary}[Remark]{Corollary}
\begin{frontmatter}
\title{Asymptotic analysis for bifurcating autoregressive processes via a martingale approach}
\runtitle{Asymptotic analysis for BAR processes}


\author{\fnms{Bernard} \snm{Bercu}\corref{}\ead[label=e1]{bernard.bercu@math.u-bordeaux1.fr}}
\address{IMB, 351 cours de la Lib\'eration, F33405 Talence, France\\ \printead{e1}\\\printead{e2}\\\printead{e3}}
\affiliation{Universit\'e de Bordeaux, IMB, CNRS, UMR 5251,\\ and INRIA Bordeaux, team CQFD, France}
\and
\author{\fnms{Beno\^\i te} \snm{de~Saporta}\ead[label=e2]{saporta@math.u-bordeaux1.fr}}
\affiliation{Universit\'e de Bordeaux, GREThA, CNRS, UMR 5113,\\ IMB, CNRS, UMR 5251, and INRIA Bordeaux, team CQFD, France}
\and
\author{\fnms{Anne} \snm{G\'egout-Petit}\ead[label=e3]{anne.petit@u-bordeaux2.fr}}
\affiliation{Universit\'e de Bordeaux, IMB, CNRS, UMR 5251,\\ and INRIA Bordeaux, team CQFD, France}

\runauthor{B. Bercu, B. de~Saporta, A. G\'egout-Petit}

\begin{abstract}\quad
We study the asymptotic behavior of the least squares estimators of the unknown parameters of general 
$p$th-order bifurcating autoregressive processes. Under very weak assumptions on the driven noise of the process, namely conditional pair-wise independence and suitable moment conditions, we establish the almost sure convergence
of our estimators together with the quadratic strong law and the central limit theorem.
All our analysis relies on non-standard asymptotic results for martingales.
\end{abstract}


\end{frontmatter}


\vspace{-5ex}
\section{Introduction}
\label{section intro}
Bifurcating autoregressive (BAR) processes are an adaptation of autoregressive (AR) processes to binary tree structured data. They were first introduced by Cowan and Staudte \cite{CoSt86} for cell lineage data, where each individual in one generation gives birth to two offspring in the next generation. Cell lineage data typically consist of observations of some quantitative characteristic of the cells over several generations of descendants from an initial cell. BAR processes take into account both inherited and environmental effects to explain the evolution of the quantitative characteristic under study.\vspace{-1ex}\\

More precisely, the original BAR process is defined as follows. The initial cell is labelled $1$, and the two offspring of cell $n$ are labelled $2n$ and $2n+1$. Denote by $X_n$ the quantitative characteristic of individual $n$. Then, the first-order BAR process is given, for all $n\geq 1$, by
\begin{equation*}
\left\{\begin{array}{lcccccl}
X_{2n} & = & a &+ &bX_n &+ &\varepsilon_{2n},\\
X_{2n+1} & = & a & +& bX_n &+ &\varepsilon_{2n+1}.\\
\end{array}\right.
\end{equation*}
The noise sequence $(\varepsilon_{2n}, \varepsilon_{2n+1})$ represents environmental effects while $a,b$ are unknown real parameters with $|b|<1$. The driven noise $(\varepsilon_{2n}, \varepsilon_{2n+1})$ was originally supposed to be independent and identically distributed with normal distribution. However,  two sister cells being in the same environment early in their lives, $\varepsilon_{2n}$ and $\varepsilon_{2n+1}$ are allowed to be correlated, inducing a correlation between sister cells distinct from the correlation inherited from their mother.\vspace{-1ex}\\

Several extensions of the model have been proposed. On the one hand, we refer the reader to Huggins and Basawa \cite{HuBa99} and Basawa and Zhou \cite{BaZh04, ZhBa05b} for statistical inference on symmetric
bifurcating processes. On the other hand, higher order processes, when not only the effects of the mother but also those of the grand-mother and higher order ancestors are taken into account, have been investigated by Huggins and Basawa \cite{HuBa99}. Recently, an asymmetric model has been introduced 
by Guyon \cite{Guy07, GBPSDT05} where only the effects of the mother are considered, but sister cells are allowed to have different conditional distributions. We can also mention a recent work
of Delmas and Marsalle \cite{DeMa08} dealing with a model of asymmetric bifurcating Markov chains 
on a Galton Watson tree instead of regular binary tree.\vspace{-1ex}\\

The purpose of this paper is to carry out a sharp analysis of the asymptotic properties of the least squares (LS) estimators of the unknown parameters of general asymmetric $p$th-order BAR processes. There are several results on statistical inference and asymptotic properties of estimators for BAR models in the literature. For maximum likelihood inference on small independent trees, see Huggins and Basawa \cite{HuBa99}. For maximum likelihood inference on a single large tree, see Huggins \cite{Hug96} for the original BAR model, Huggins and Basawa \cite{HuBa00} for higher order Gaussian BAR models, and Zhou and Basawa \cite{ZhBa05b} for exponential first-order BAR processes. We also refer the reader to Zhou and Basawa \cite{ZhBa05a} for the LS parameter estimation, and to Hwang, Basawa and Yeo \cite{HwBaYeo09} for the local asymptotic
normality for BAR processes and related asymptotic inference. In all those papers, the process is supposed to be stationary. Consequently, $X_n$ has a time-series representation involving an holomorphic function. In Guyon \cite{Guy07}, the LS estimator is also investigated, but the process is not stationary, and the author makes intensive use of the tree structure and Markov chain theory. Our goal is to improve and extend the previous results of Guyon \cite{Guy07} via a martingale approach. As previously done by Basawa and Zhou \cite{BaZh04, ZhBa05a, ZhBa05b} we shall make use of the strong law of large numbers \cite{Duflo97} as well as the central limit theorem \cite{HaHe80, Ham94} for martingales. It will allow us to go further in the analysis of general $p$th-order BAR processes. We shall establish the almost sure convergence of the LS estimators together with the quadratic strong law and the central limit theorem.\vspace{-1ex}\\

The paper is organised as follows. Section~\ref{section model} is devoted to the presentation of the asymmetric $p$th-order BAR process under study, while Section~\ref{section LS} deals with the LS estimators of the unknown parameters. In Section~\ref{sectionmartingale}, we explain our strategy based on martingale theory. Our main results about the asymptotic properties of the LS estimators are given in Section~\ref{section results}. More precisely, we shall establish the almost sure convergence, the quadratic strong law (QSL) and the central limit theorem (CLT) for the LS estimators. The proof of our main
results are detailed in Sections 6 to 10, the more technical ones being gathered in the appendices.

\section{Bifurcating autoregressive processes}
\label{section model}
In all the sequel, let $p$ be a non-zero integer. We consider the asymmetric BAR($p$) process given, for all $n\geq 2^{p-1}$, by
\begin{equation}\label{defbar}
\left\{
\begin{array}{lccccccl}
X_{2n} & = & a_{0} &+&\sum_{k=1}^p a_k X_{[\frac{n}{2^{k-1}}]} &+&\varepsilon_{2n}, \\
X_{2n+1} & = & b_{0} &+&\sum_{k=1}^p b_k X_{[\frac{n}{2^{k-1}}]} &+&\varepsilon_{2n+1},
\end{array}\right.
\end{equation}
where $[x]$ stands for the largest integer less than or equal to $x$. 
The initial states $\{X_k,\ 1\leq k\leq 2^{p-1}-1\}$ are the ancestors while $(\varepsilon_{2n},\varepsilon_{2n+1})$ is the driven noise of the process. The parameters 
$(a_0,a_1,\ldots a_p)$ and $(b_0,b_1,\ldots,b_p)$ are unknown real numbers. 
The BAR($p$) process can be rewritten in 
the abbreviated vector form given, for all $n\geq 2^{p-1}$, by
\begin{equation}\label{defbarmatrix}
\left\{
\begin{array}{lcccl}
\dX_{2n} & = & A\dX_n&+&\eta_{2n}, \\
\dX_{2n+1} & = & B\dX_n&+&\eta_{2n+1},
\end{array}\right.
\end{equation}
where the regression vector $\dX_n=(X_n, X_{[\frac{n}{2}]},\ldots,X_{[\frac{n}{2^{p-1}}]})^t$,
$\eta_{2n}=(a_0+\varepsilon_{2n})e_1$, $\eta_{2n+1}=(b_0+\varepsilon_{2n+1})e_1$ with $e_1=(1,0,\ldots,0)^t\in\dR^p$. Moreover, $A$ and $B$ are the $p\times p$ companion matrices
\begin{equation*}
A=\left(\begin{array}{cccc}
a_1&a_2&\cdots&a_p\\
1&0&\cdots&0\\
0&\ddots&\ddots&\vdots\\
0&0&1&0
\end{array}
\right),\qquad B=\left(\begin{array}{cccc}
b_1&b_2&\cdots&a_p\\
1&0&\cdots&0\\
0&\ddots&\ddots&\vdots\\
0&0&1&0
\end{array}
\right).
\end{equation*}
This process is a direct generalization of the symmetric BAR($p$) process studied by
Huggins, Basawa and Zhou \cite{HuBa99, ZhBa05a}.
One can also observe that, in the particular case $p=1$, it is 
the asymmetric BAR process studied by Guyon
\cite{Guy07, GBPSDT05}. In all the sequel, we shall assume that $\mathbb{E}[X_k^8]<\infty$ 
for all $1\leq k\leq 2^{p-1}-1$ and that matrices $A$ and $B$ satisfy the contracting property
\begin{equation*}
\beta=\max\{\|A\|, \|B\|\}<1,
\end{equation*}
where $\|A\|=\sup\{\|Au\|, \ u\in \dR^p \text{ with } \|u\|=1\}$. 

\begin{figure}[h]
\centering
\includegraphics[height=8cm, angle=-90]{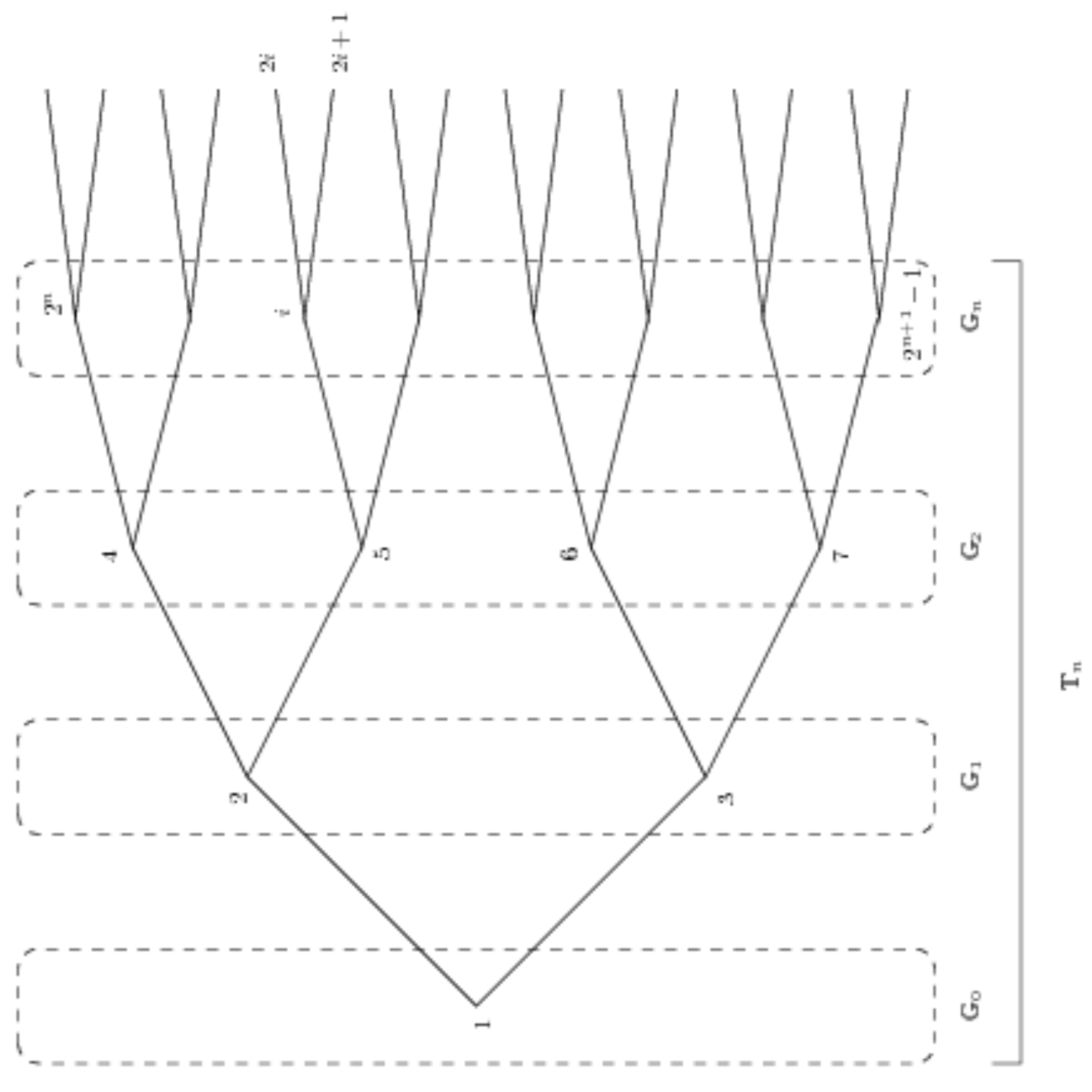}
\caption{The tree associated with the bifurcating auto-regressive process.\label{arbre}}
\end{figure}

As explained in the introduction, one can see this BAR($p$) process as a $p$th-order 
autoregressive process on a binary tree,
where each vertex represents an individual or cell, vertex $1$ being the original ancestor, see Figure \ref{arbre} for an illustration. For all $n\geq 1$, denote the $n$th generation by
$$\dG_n = \{2^n, 2^n+1,\ldots,2^{n+1}-1\}.$$
In particular, $\dG_0 = \{1\}$ is the initial generation and $\dG_1 = \{2,3\}$ is the first generation of offspring from the first ancestor. Let $\dG_{r_n}$ be the generation of individual $n$, which means that $r_n=\log_2(n)$. Recall that the two offspring of individual $n$ are labelled $2n$ and $2n+1$, or conversely, the mother of individual $n$ is $[n/2]$. More generally, the ancestors of individual $n$ are $[n/2], [n/2^2],\ldots,\ [n/2^{r_n}]$. Furthermore, denote by
$$\dT_n = \bigcup_{k=0}^n\dG_k$$
the sub-tree of all individuals from the original individual up to the $n$th generation. It is clear 
that the cardinality $|\dG_n|$ of $\dG_n$ is $2^n$ while that of $\dT_n$ is $|\dT_n|=2^{n+1}-1$.
Finally, we denote by $\dT_{n,p}=\{ k \in \dT_n , k\geq 2^{p} \}$
the sub-tree of all individuals up to the $n$th generation without $\dT_{p-1}$.
One can observe that, for all $n\geq 1$, $\dT_{n,0}= \dT_{n}$ and, for all $p\geq 1$, 
$\dT_{p,p}=\dG_{p}$.  

\section{Least-squares estimation}
\label{section LS}
The BAR($p$) process (\ref{defbar}) can be rewritten, for all $n\geq 2^{p-1}$, in the matrix form
\begin{equation}\label{defbarmat}
Z_n=\theta^tY_n+V_n
\end{equation}
where
\begin{equation*}
Z_n=\left(
    \begin{array}{l}
     X_{2n}\\
     X_{2n+1}
    \end{array}\right),
\hspace{0.5cm}
Y_n=\left(
    \begin{array}{l}
     1\\
     \dX_{n}
    \end{array}\right),
\hspace{0.5cm}
V_n=\left(
     \begin{array}{l}
     \varepsilon_{2n}\\
     \varepsilon_{2n+1}
    \end{array}\right),
\end{equation*}
and the $(p+1)\times2$ matrix parameter $\theta$ is given by
\begin{equation*}
\theta=
\left(\begin{array}{cc}
a_0 & b_0  \\
a_1 & b_1\\
\vdots&\vdots\\
a_p&b_p
\end{array}\right).
\end{equation*}
Our goal is to estimate $\theta$ from the observation of all individuals up to the $n$th generation that is the complete sub-tree $\dT_n$. Each new generation $\dG_n$ contains half the global available information.
Consequently, we shall show that observing the whole tree
$\dT_n$ or only generation $\dG_n$ is almost the same.
We propose to make use of the standard 
LS estimator $\wh{\theta}_n$ which minimizes
$$
\Delta_n(\theta)=\frac{1}{2}\sum_{k\in \dT_{n-1,p-1}}\parallel Z_k -\theta^t Y_k \parallel^2.
$$
Consequently, we obviously have for all $n\geq p$
\begin{equation}
\label{defLS}
\wh{\theta}_n=S_{n-1}^{-1}\sum_{k\in \dT_{n-1,p-1}}Y_kZ_k^t,
\end{equation}
where the $(p+1)\times(p+1)$ matrix $S_n$ is defined as
\begin{equation*}
S_n =\sum_{k \in \dT_{n,p-1}}Y_kY_k^t =
\sum_{k \in \dT_{n,p-1}}
\left(\begin{array}{cc}
1 & \dX_k^t \\
\dX_k  & \dX_k\dX_k^t
\end{array}\right).
\end{equation*}
In the special case where $p=1$, $S_n$ simply reduces to
\begin{equation*}
S_n =
\sum_{k \in \dT_n}
\left(\begin{array}{cc}
1 & X_k \\
X_k  & X_k^2
\end{array}\right).
\end{equation*}
In order to avoid useless invertibility assumption, we shall assume, without loss of generality, that for all $n\geq {p-1}$, $S_n$ is invertible. Otherwise, we only have to add the identity matrix $\rI_{p+1}$ to $S_n$. In all what follows, we shall make a slight abuse of notation by identifying $\theta$ as well as $\wh{\theta}_n$ to
$$
\text{vec}(\theta)=
\left(
    \begin{array}{c}
a_0\\
\vdots\\
a_p\\
b_0\\
\vdots\\
b_p
\end{array}\right)
\hspace{1cm}\text{and}\hspace{1cm}
\text{vec}(\wh{\theta}_n)=
\left(
    \begin{array}{c}
     \wh{a}_{0,n}\\
\vdots\\
     \wh{a}_{p,n}\\
     \wh{b}_{0,n}\\
\vdots\\
     \wh{b}_{p,n}
    \end{array}\right).
$$
The reason for this change will be explained in Section~\ref{sectionmartingale}. Hence, we readily deduce from (\ref{defLS}) that
\begin{eqnarray*}
\wh{\theta}_n&=&
(\rI_2 \otimes S^{-1}_{n-1}) \sum_{k\in \dT_{n-1,p-1}}\text{vec}\left(Y_kZ_k^t\right)\nonumber\\
&=&(\rI_2 \otimes S^{-1}_{n-1})
\sum_{k \in \mathbb{T}_{n-1,p-1}}
\left( \begin{array}{c}
X_{2k}  \\
X_k\dX_{2k} \\
X_{2k+1} \\
X_k\dX_{2k+1}
\end{array}\right),
\end{eqnarray*}
where $\otimes$ stands for the matrix Kronecker product. Consequently, it follows from (\ref{defbarmat}) 
that
\begin{eqnarray} \label{diff}
\wh{\theta}_n-\theta &=&(\rI_2 \otimes S^{-1}_{n-1})  \sum_{k\in \dT_{n-1,p-1}}\text{vec}\left(Y_kV_k^t\right)\nonumber\\
&=&(\rI_2 \otimes S^{-1}_{n-1})
\sum_{k \in \dT_{n-1,p-1}}
\left( \begin{array}{c}
\veps_{2k}  \\
\veps_{2k}\dX_k \\
\veps_{2k+1} \\
\veps_{2k+1}\dX_k
\end{array}\right).
\end{eqnarray}
Denote by $\dF=(\cF_n)$ the natural filtration associated with the BAR($p$) process, which means that $\cF_n$ is the $\sigma$-algebra generated by all individuals up to the $n$th generation, 
$\cF_n = \sigma\{X_k, k\in \dT_n\}$. In all the sequel, we shall make use of the five following moment hypotheses.
\begin{description}
\item[(H.1)] One can find $\sigma^2>0$ such that, for all $n\geq {p-1}$ and for all $k\in \dG_{n+1}$, 
$\veps_k$ belongs to $L^2$ with
\begin{equation*}
\dE[\veps_k|\cF_n] = 0
\hspace{1cm}\text{and}\hspace{1cm}
\dE[\veps_k^2|\cF_n]=\sigma^2
\hspace{1cm}\text{a.s.}
\end{equation*}
\item[(H.2)] It exists $|\rho|<\sigma^2$ such that, for all $n\geq {p-1}$ and for all different $k,l\in \dG_{n+1}$ with $[k/2]=[l/2]$, 
\begin{equation*}
\dE[\veps_k\veps_l|\cF_n] = \rho \hspace{1cm}\text{a.s.}
\end{equation*}
Otherwise, $\veps_k$ and $\veps_l$ are conditionally independent given $\cF_n$.
\item[(H.3)] For all $n\geq {p-1}$ and for all $k\in \dG_{n+1}$,
$\veps_k$ belongs to $L^4$ and
\begin{equation*}
\sup_{n\geq {p-1}}\sup_{k\in \dG_{n+1}}
\dE[\veps_k^4|\cF_n]<\infty
\hspace{1cm}\text{a.s.}
\end{equation*}
\item[(H.4)] One can find $\tau^4>0$ such that, for all $n\geq {p-1}$ and 
for all $k\in \dG_{n+1}$, 
\begin{equation*}
\dE[\veps_k^4|\cF_n]=\tau^4
\hspace{1cm}\text{a.s.}
\end{equation*}
and, for $\nu^2<\tau^4$ and for all different $k,l\in \dG_{n+1}$ with $[k/2]=[l/2]$ 
\begin{equation*}
\dE[\veps_{2k}^2\veps_{2k+1}^2|\cF_n] = \nu^2 \hspace{1cm}\text{a.s.}
\end{equation*}
\item[(H.5)] For all $n\geq {p-1}$ and for all $k\in \dG_{n+1}$,
$\veps_k$ belongs to $L^8$ with
\begin{equation*}
\sup_{n\geq {p-1}}\sup_{k\in \dG_{n+1}}
\dE[\veps_k^8|\cF_n]<\infty
\hspace{1cm}\text{a.s.}
\end{equation*}
\end{description}
\begin{Remark}$\pt$
In contrast with \cite{ZhBa05a}, one can observe that we do not assume that $(\veps_{2n},\veps_{2n+1})$ is a sequence of independent and identically distributed bi-variate random vectors. The price to pay for giving up this iid assumption is higher moments, namely assumptions \emph{(\textbf{H.3})} and \emph{(\textbf{H.5})}. Indeed we need them to make use of the strong law of large numbers and 
the central limit theorem for martingales. However, we do not require
any normality assumption on $(\veps_{2n},\veps_{2n+1})$. Consequently, our assumptions are
much weaker than the existing ones in previous literature.
\end{Remark}
We now turn to the estimation of the parameters $\sigma^2$ and $\rho$. On the one hand, we propose to estimate the conditional variance $\sigma^2$ by
\begin{equation}
\label{sigma2}
\wh{\sigma}^2_n = \frac{1}{2 |\dT_{n-1}|}\sum_{k \in \dT_{n-1,p-1}} \parallel \wh{V}_k \parallel^2
=\frac{1}{2 |\dT_{n-1}|}\sum_{k \in\dT_{n-1,p-1}}(\wh{\veps}_{2k}^2 +   \wh{\varepsilon}_{2k+1}^2)
\end{equation}
where for all $n\geq {p-1}$ and for all $k\in\mathbb{G}_n$, $\wh{V}_k^t=(\wh{\veps}_{2k}, \wh{\veps}_{2k+1})$ with
\begin{equation*}
\left\{
\begin{array}{lccccccl}
\wh{\veps}_{2k} & = & X_{2k} &-& \wh{a}_{0,n} &-& \sum_{i=1}^p \wh{a}_{i,n}X_{[\frac{k}{2^{i-1}}]},\\
\wh{\veps}_{2k+1} & = & X_{2k+1} &-& \wh{b}_{0,n} &-& \sum_{i=1}^p \wh{b}_{i,n}X_{[\frac{k}{2^{i-1}}]}.
\end{array}\right.
\end{equation*}
One can observe that, on the above equations, we make use of only the past observations for the
estimation of the parameters. This will be crucial in the asymptotic analysis. 
On the other hand, we estimate the conditional covariance $\rho$ by
\begin{equation}
\label{rho}
\wh{\rho}_n = \frac{1}{|\dT_{n-1}|}\sum_{k \in\dT_{n-1,p-1}}
\wh{\veps}_{2k} \wh{\veps}_{2k+1}.
\end{equation}

\section{Martingale approach}
\label{sectionmartingale}
In order to establish all the asymptotic properties of our estimators, we shall make use of a martingale approach. It allows us
to impose a very smooth restriction on the driven noise $(\veps_n)$ compared with the previous results in the literature. As a matter of fact, we only assume suitable moment conditions on $(\veps_n)$ and that $(\varepsilon_{2n},\varepsilon_{2n+1})$ are
conditionally independent, while it is assumed in \cite{ZhBa05a} that $(\varepsilon_{2n},\varepsilon_{2n+1})$ is a sequence of independent identically distributed random vectors. For all $n\geq p$, denote
\begin{equation*}
M_n= \sum_{k \in \dT_{n-1,p-1}} \left( \begin{array}{c}
\veps_{2k}  \\
\veps_{2k}\dX_k \\
\veps_{2k+1} \\
\veps_{2k+1}\dX_k
\end{array}\right)\in\dR^{2(p+1)}.
\end{equation*}
Let $\Sigma_n=\rI_2 \otimes S_{n}$, and note that 
$\Sigma^{-1}_{n}=\rI_2 \otimes S_{n}^{-1}$. For all $n\geq p$, we can
thus rewrite (\ref{diff}) as
\begin{equation}
\label{thetadiff} \widehat{\theta}_n-\theta  = \Sigma^{-1}_{n-1}
M_{n}.
\end{equation}
The key point of our approach is that $(M_n)$ is a martingale. Most of all the asymptotic results for martingales were established for vector-valued martingales. That is the reason why we have chosen to make use of vector notation in  Section~\ref{section LS}. In order to show that $(M_n)$ is a martingale adapted to the filtration $\dF=(\cF_n)$, we rewrite it in a compact form. Let $\Psi_n=\rI_2\otimes \Phi_n$, where $\Phi_n$ is the rectangular matrix of dimension $(p+1)\times \delta_n$, with $\delta_n=2^n$, given by
\begin{equation*}
\Phi_n= \left(\begin{array}{cccc}
1 & 1 & \cdots & 1 \vspace{1ex}\\
\dX_{2^n} & \dX_{2^n + 1} & \cdots & \dX_{2^{n+1}-1}
\end{array}\right).
\end{equation*}
It contains the individuals of generations $\dG_{n-p+1}$ up to $\dG_{n}$ and is also the collection of all $Y_k$, $k\in\mathbb{G}_n$. Let $\xi_n$ be the random vector of dimension $\delta_n$
\begin{equation*}
\xi_n= \left( \begin{array}{cccccccc}
\veps_{2^n}  \\
\veps_{2^n+2} \\
\vdots \\
\veps_{2^{n+1}-2} \\
\veps_{2^n+1}  \\
\veps_{2^n+3} \\
\vdots \\
\veps_{2^{n+1}-1} \\
\end{array}\right).
\vspace{2ex}
\end{equation*}
The vector $\xi_n$ gathers the noise variables of generation $\dG_n$. The special ordering separating odd and even indices is tailor-made so that $M_n$ can be written as
\begin{equation*}
M_n=\sum_{k=p}^n\Psi_{k-1}\xi_k.
\end{equation*}
By the same token, one can observe that
\begin{equation*}
S_n =\sum_{k=p-1}^n\Phi_k\Phi_k^t
\hspace{1cm}\text{and}\hspace{1cm}
\Sigma_n =\sum_{k=p-1}^n\Psi_k\Psi_k^t.
\vspace{2ex}
\end{equation*}
Under (\textbf{H.1}) and (\textbf{H.2}), we clearly have for all $n\geq 0$, $\dE[\xi_{n+1}|\cF_n]=0$ and $\Psi_n$ is
$\cF_n$-measurable. In addition, it is not hard to see that for all $n\geq 0$, $\dE[\xi_{n+1}\xi_{n+1}^t|\cF_n]=\Gamma \otimes
\rI_{\delta_n}$ where $\Gamma$ is the covariance matrix associated with $(\veps_{2n},\veps_{2n+1})$
\begin{equation*}
\Gamma=
\left(\begin{array}{cc}
\sigma^2 & \rho \vspace{1ex}\\
\rho & \sigma^2
\end{array}\right).
\end{equation*}
We shall also prove that $(M_n)$ is a square integrable martingale. 
Its increasing process is given for all $n\geq p+1$ by
\begin{equation*}
\cM_n=\sum_{k=p-1}^{n-1}\Psi_{k}(\Gamma \otimes\rI_{\delta_k})\Psi_{k}^t
=\Gamma \otimes\sum_{k=p-1}^{n-1}\Phi_{k}\Phi_k^t
= \Gamma \otimes S_{n-1}.
\end{equation*}
It is necessary to establish the convergence of $S_n$, properly normalized, in order to prove the asymptotic results for the BAR($p$) estimators $\wh{\theta}_n$, $\wh{\sigma}^2_n$ and $\wh{\rho}_n$. One can observe that the sizes of $\Psi_{n}$ and $\xi_n$ are not fixed and double at each generation. This is why we have to adapt the proof of vector-valued martingale convergence given in \cite{Duflo97} to our framework.

\section{Main results}
\label{section results}
We now state our main results, first on the martingale $(M_n)$ and then on our estimators.
\begin{Proposition}$\pt$
\label{mainlemma}
Assume that $(\veps_n)$ satisfies \emph{(\textbf{H.1})} to \emph{(\textbf{H.3})}.
Then, we have
\begin{equation}
\label{cvgSn}
\lim_{n\rightarrow \infty} \frac{S_{n}}{|\dT_n|} = L \hspace{1cm} \text{a.s.}
\end{equation}
where $L$ is a positive definite matrix specified in Section~\ref{section Sn}.
\end{Proposition}
This result is the keystone of our asymptotic analysis. It enables us to prove sharp asymptotic properties for $(M_n)$.
\begin{Theorem}\label{th Mn}
Assume that $(\veps_n)$ satisfies \emph{(\textbf{H.1})} to \emph{(\textbf{H.3})}. Then, we have
\begin{equation}\label{th Mn1}
M_n^t\Sigma_{n-1}^{-1}M_n=\cO(n)\hspace{1cm} \text{a.s.}
\end{equation}
In addition, we also have
\begin{equation}\label{th Mn3}
\lim_{n\rightarrow\infty}\frac{1}{n}\sum_{k=p}^n M_k^{t}\Sigma_{k-1}^{-1}M_k=2(p+1)\sigma^2\hspace{1cm} \text{a.s.}
\end{equation}
Moreover, if $(\veps_n)$ satisfies \emph{(\textbf{H.4})} and \emph{(\textbf{H.5})}, 
we have the central limit theorem
\begin{equation}\label{th Mn4}
\frac{1}{\sqrt{|\dT_{n-1}|}}M_n\build{\longrightarrow}_{}^{{\mbox{\calcal L}}}
\cN(0,\Gamma\otimes L).
\end{equation}
\end{Theorem}
From the asymptotic properties of $(M_n)$, we deduce the asymptotic behavior of our estimators. Our first result deals with the almost sure asymptotic properties of the LS estimator $\wh{\theta}_n$.
\begin{Theorem}$\pt$
\label{thmaptheta}
Assume that $(\veps_n)$ satisfies \emph{(\textbf{H.1})} to \emph{(\textbf{H.3})}.
Then, $\wh{\theta}_{n}$ converges almost surely to $\theta$
with the rate of convergence
\begin{equation}
\label{aptheta1}
\parallel \widehat{\theta}_{n}-\theta \parallel^{2}=
\cO \left(\frac{\log |\dT_{n-1}|}{|\dT_{n-1}|} \right)
\hspace{1cm}\text{a.s.}
\end{equation}
In addition, we also have the quadratic strong law
\begin{equation}\label{aptheta2}
\lim_{n\rightarrow \infty} \frac{1}{n}\sum_{k=1}^n |\dT_{k-1}|
(\widehat{\theta}_{k}-\theta)^t\Lambda (\widehat{\theta}_{k}-\theta)= 2(p+1)\sigma^2 \hspace{1cm} \text{a.s.}
\end{equation}
where $\Lambda=\rI_2\otimes L$.
\end{Theorem}

\noindent Our second result is devoted to the almost sure asymptotic properties of
the variance and covariance estimators $\wh{\sigma}^2_n$ and $\wh{\rho}_n$. Let
\begin{equation*}
\sigma^2_n =\frac{1}{2 |\dT_{n-1}|}\sum_{k \in\dT_{n-1,p}}({\veps}_{2k}^2 + {\varepsilon}_{2k+1}^2)
\hspace{0.2cm}\text{and}\hspace{0.2cm}
{\rho}_n = \frac{1}{|\dT_{n-1}|}\sum_{k \in\dT_{n-1,p}} {\veps}_{2k}{\veps}_{2k+1}.
\end{equation*}

\begin{Theorem}$\pt$
\label{thmapsigmarho}
Assume that $(\veps_n)$ satisfies \emph{(\textbf{H.1})} to \emph{(\textbf{H.3})}.
Then, $\wh{\sigma}^2_n$ converges almost surely to $\sigma^2$. More precisely,
\begin{equation}\label{apsigma1}
\lim_{n\rightarrow\infty}\frac{1}{n}\sum_{k\in\dT_{n-1,p}}(\wh{\varepsilon}_{2k}-\varepsilon_{2k})^2+(\wh{\varepsilon}_{2k+1}-\varepsilon_{2k+1})^2=2(p+1)\sigma^2
\hspace{1cm}\text{a.s.}
\end{equation}
\begin{equation}
\label{apsigma2}
\lim_{n\rightarrow\infty}\frac{|\dT_{n}|}{n}
(\wh{\sigma}^2_n-{\sigma}^2_n)=2(p+1)\sigma^2
\hspace{1cm}\text{a.s.}
\end{equation}
In addition, $\wh{\rho}_n$ converges almost surely to $\rho$
\begin{equation}
\label{aprho1}
\lim_{n\rightarrow\infty}\frac{1}{n}\sum_{k\in\dT_{n-1,p}}(\wh{\varepsilon}_{2k}-\varepsilon_{2k})(\wh{\varepsilon}_{2k+1}-\varepsilon_{2k+1})=(p+1)\rho
\hspace{1cm}\text{a.s.}
\end{equation}
\begin{equation}
\label{aprho2}
\lim_{n\rightarrow\infty}\frac{|\dT_{n}|}{n}
(\wh{\rho}_n-{\rho}_n)=2(p+1)\rho
\hspace{1cm}\text{a.s.}
\end{equation}
\end{Theorem}

\noindent
Our third result concerns the asymptotic normality for all our estimators
$\wh{\theta}_n$, $\wh{\sigma}^2_n$ and $\wh{\rho}_n$.

\begin{Theorem}$\pt$
\label{thmCLT}
Assume that $(\veps_n)$ satisfies \emph{(\textbf{H.1})} to \emph{(\textbf{H.5})}.
Then, we have the central limit theorem
\begin{equation}
\label{CLTtheta}
\sqrt{|\dT_{n-1}|} (\wh{\theta}_{n}-\theta)
\build{\longrightarrow}_{}^{{\mbox{\calcal L}}}
\cN(0,\Gamma\otimes L^{-1}).
\end{equation}
In addition, we also have
\begin{equation}
\label{CLTsigma}
\sqrt{|\dT_{n-1}|} (\wh{\sigma}^2_n-{\sigma}^2)
\build{\longrightarrow}_{}^{{\mbox{\calcal L}}}
\cN\Bigl(0,\frac{\tau^4-2\sigma^4+\nu^2}{2}\Bigr)
\end{equation}
and
\begin{equation}
\label{CLTrho}
\sqrt{|\dT_{n-1}|} (\wh{\rho}_n-{\rho})
\build{\longrightarrow}_{}^{{\mbox{\calcal L}}}
\cN(0,\nu^2-\rho^2).
\end{equation}
\end{Theorem}
The rest of the paper is dedicated to the proof of our main results. We start by giving laws of large numbers for the noise sequence $(\varepsilon_n)$ in Section~\ref{sectionLLnEps}. In Section~\ref{section Sn}, we give the proof of
Proposition~\ref{mainlemma}. Sections~\ref{section th1},~\ref{section th2} and~\ref{section th3} are devoted to the proofs
of Theorems~\ref{thmaptheta},~\ref{thmapsigmarho} and~\ref{thmCLT}, respectively. The more technical proofs, including that of Theorem~\ref{th Mn}, are postponed to the Appendices.

\section{Laws of large numbers for the noise sequence}
\label{sectionLLnEps}
We first need to establish strong laws of large numbers for the noise sequence $(\varepsilon_n)$. These results will be useful in all the sequel. We will extensively use the strong law of large numbers for locally square integrable real martingales given in Theorem 1.3.15 of \cite{Duflo97}.

\begin{Lemma}$\pt$
\label{lemLLNeps}
Assume that $(\veps_n)$ satisfies \emph{(\textbf{H.1})} and \emph{(\textbf{H.2})}.
Then
\begin{equation}
\label{LLNeps1}
\lim_{n \rightarrow + \infty}\frac{1}{|\dT_n|}\sum_{k \in \dT_{n,p}}\veps_k=0
\hspace{1cm}\text{a.s.}
\end{equation}
In addition, if \emph{(\textbf{H.3})} holds, we also have
\begin{equation}
\label{LLNeps2}
\lim_{n \rightarrow + \infty}\frac{1}{|\dT_n|}\sum_{k \in\dT_{n,p}}\veps_k^2=\sigma^2
\hspace{1cm}\text{a.s.}
\end{equation}
and
\begin{equation}
\label{LLNeps3}
\lim_{n \rightarrow + \infty}\frac{1}{|\dT_{n-1}|}\sum_{k \in\dT_{n-1,p-1}}\veps_{2k}\veps_{2k+1}=\rho
\hspace{1cm}\text{a.s.}
\end{equation}
\end{Lemma}

\noindent\textbf{Proof :} On the one hand, let
\begin{equation*}
P_n=\sum_{k \in \dT_{n,p}}\veps_k=\sum_{k=p}^n\sum_{i \in \dG_k} \veps_i.
\end{equation*}
We have
\begin{equation*}
\Delta P_{n+1}=P_{n+1}-P_n=\sum_{k \in \dG_{n+1}}\veps_k.
\end{equation*}
Hence, it follows from (\textbf{H.1}) and (\textbf{H.2})
that $(P_n)$ is a square integrable real martingale with increasing process
\begin{equation*}
<\!P\!>_n=(\sigma^2+\rho)\sum_{k=p}^n|\dG_{k}|=(\sigma^2+\rho)(|\dT_{n}|-|\dT_{p-1}|).
\end{equation*}
Consequently, we deduce from Theorem 1.3.15 of \cite{Duflo97} that $P_n=o(<\!P\!>_n)$ a.s. which implies (\ref{LLNeps1}).
On the other hand, denote
\begin{equation*}
Q_n=\sum_{k=p}^n \frac{1}{|\dG_k|}\sum_{i \in \dG_k}e_i,
\end{equation*}
where $e_n=\veps_n^2 - \sigma^2$. We have
\begin{equation*}
\Delta Q_{n+1}=Q_{n+1}-Q_n=\frac{1}{|\dG_{n+1}|}\sum_{k \in \dG_{n+1}}e_k.
\end{equation*}
First of all, it follows from (\textbf{H.1}) that for all $k \in \dG_{n+1}$, $\dE[e_k|\cF_n] = 0$ a.s. In addition,
for all different $k,l\in \dG_{n+1}$ with $[k/2]\neq[l/2]$,
\begin{equation*}
\label{majnoisecove}
\dE[e_ke_l|\cF_n] =0
\hspace{1cm}\text{a.s.}
\end{equation*}
thanks to the conditional independence given by (\textbf{H.2}). Furthermore, we readily deduce from (\textbf{H.3}) that
\begin{equation*}
\label{majnoisemom2e}
\sup_{n\geq p-1}\sup_{k\in \dG_{n+1}}
\dE[e_k^2|\cF_n]<\infty
\hspace{1cm}\text{a.s.}
\end{equation*}
Therefore, $(Q_n)$ is a square integrable real martingale with increasing process
\begin{eqnarray*}
<\!Q\!>_n
& \! \leq \! & 2 \sup_{p-1\leq k \leq n-1}\sup_{i\in \dG_{k+1}}
\dE[e_i^2|\cF_k]\sum_{j=p}^n \frac{1}{|\dG_{j}|} \hspace{0.5cm}\text{a.s.} \\
& \! \leq \! & 2 \sup_{p-1\leq k \leq n-1}\sup_{i\in \dG_{k+1}}
\dE[e_i^2|\cF_k]\sum_{j=p}^n \Bigl(\frac{1}{2}\Bigr)^j \hspace{0.5cm}\text{a.s.} \\
& & \\
& \! \leq \! & 2 \sup_{p-1\leq k \leq n-1}\sup_{i\in \dG_{k+1}}
\dE[e_i^2|\cF_k] <\infty
\hspace{1cm}\text{a.s.}
\end{eqnarray*}
Consequently, we obtain from the strong law of large numbers for martingales
that $(Q_n)$ converges almost surely. Finally, as $(|\dG_n|)$ is a positive real sequence which increases to infinity, we find from Lemma~\ref{lemKro} in Appendix~\ref{appendixA} that
\begin{equation*}
\sum_{k=p}^n \sum_{i \in \dG_k}e_i=o(|\dG_n|)\hspace{1cm}\text{a.s.}
\end{equation*}
leading to
\begin{equation*}
\sum_{k=p}^n \sum_{i \in \dG_k}e_i=o(|\dT_n|)\hspace{1cm}\text{a.s.}
\end{equation*}
as $|\dT_n|-1=2|\dG_n|$, which implies (\ref{LLNeps2}).
We also establish (\ref{LLNeps3}) in a similar way. As a matter of fact, let
\begin{equation*}
R_n=\sum_{k=p}^n \frac{1}{|\dG_{k-1}|}\sum_{i \in \dG_{k-1}}(\veps_{2i}\veps_{2i+1}-\rho).
\end{equation*}
Then, $(R_n)$ is a square integrable real martingale which converges almost surely, leading to (\ref{LLNeps3}).
\hspace{\stretch{1}}$\Box$\\

\begin{Remark}$\pt$\label{rqCVeps2k}
Note that via Lemma \ref{lem tree-generation}
\begin{eqnarray*}
\lim_{n \rightarrow + \infty}\frac{1}{|\dG_n|}\sum_{k \in\dG_n}\veps_{2k}=0,&&
\lim_{n \rightarrow + \infty}\frac{1}{|\dG_n|}\sum_{k \in\dG_n}\veps_{2k+1}=0
\hspace{1cm}\text{a.s.}\\
\lim_{n \rightarrow + \infty}\frac{1}{|\dG_n|}\sum_{k \in\dG_n}\veps_{2k}^2=\sigma^2,&&
\lim_{n \rightarrow + \infty}\frac{1}{|\dG_n|}\sum_{k \in\dG_n}\veps_{2k+1}^2=\sigma^2
\hspace{1cm}\text{a.s.}
\end{eqnarray*}
In fact, each new generation contains half the global available information, observing the whole tree
$\dT_n$ or only generation $\dG_n$ is essentially the same.
\end{Remark}

For the CLT, we will also need the convergence of higher moments of the driven noise $(\veps_n)$.
\begin{Lemma}$\pt$
\label{lemLLEps4TCL}
Assume that $(\veps_n)$ satisfies \emph{(\textbf{H.1})} to \emph{(\textbf{H.5})}.
Then, we have
\begin{equation*}
\lim_{n \rightarrow + \infty}\frac{1}{|\dT_n|}\sum_{k \in\dT_{n,p}}\veps_k^4=\tau^4
\hspace{1cm}\text{a.s.}
\end{equation*}
and
\begin{equation*}
\lim_{n \rightarrow + \infty}\frac{1}{|\dT_{n-1}|}\sum_{k \in\dT_{n-1,p-1}}
\veps_{2k}^2\veps_{2k+1}^2=\nu^2
\hspace{1cm}\text{a.s.}
\end{equation*}
\end{Lemma}

\noindent\textbf{Proof :} The proof is left to the reader as it follows essentially the same lines as the proof of Lemma~\ref{lemLLNeps} using the square integrable real martingales
\begin{equation*}
Q_n=\sum_{k=p}^n \frac{1}{|\dG_k|}\sum_{i \in \dG_k}(\veps_i^4-\tau^4)
\end{equation*}
and
\begin{equation*}
R_n=\sum_{k=p}^n \frac{1}{|\dG_{k-1}|}\sum_{i \in \dG_{k-1}}(\veps_{2i}^2\veps_{2i+1}^2-\nu^2).
\end{equation*}

\begin{Remark}$\pt$\label{rqCVeps4}
Note that again via Lemma \ref{lem tree-generation} 
\begin{eqnarray*}
\lim_{n \rightarrow + \infty}\frac{1}{|\dG_n|}\sum_{k \in\dG_n}\veps_{2k}^4=\tau^4&
\textrm{\emph{and}}&
\lim_{n \rightarrow + \infty}\frac{1}{|\dG_n|}\sum_{k \in\dG_n}\veps_{2k+1}^4=\tau^4
\hspace{1cm}\text{a.s.}
\end{eqnarray*}
\end{Remark}

\section{Proof of Proposition~\ref{mainlemma}}
\label{section Sn}

Proposition~\ref{mainlemma} is a direct application of the two following lemmas
which provide two strong laws of large numbers for the sequence of random vectors $(\dX_n)$.

\begin{Lemma}$\pt$
\label{lemlimsumX}
Assume that $(\veps_n)$ satisfies \emph{(\textbf{H.1})} and \emph{(\textbf{H.2})}.
Then, we have
\begin{equation}
\label{cvgsumX}
\lim_{n \rightarrow + \infty}\frac{1}{|\dT_n|} \sum_{k \in\dT_{n,p}}\dX_k = \lambda = \overline{a}(I_{p}-\overline{A})^{-1}e_1
\hspace{1cm}\text{a.s.}
\end{equation}
where $\overline{a}=(a_0+b_0)/2$ and $\overline{A}$ is the 
mean of the companion matrices
$$\overline{A}=\frac{1}{2}(A+B).$$
\end{Lemma}

\begin{Lemma}$\pt$
\label{lemlimsumXX}
Assume that $(\veps_n)$ satisfies \emph{(\textbf{H.1})} to \emph{(\textbf{H.3})}. Then, we have
\begin{equation}
\label{cvgsumXX}
\lim_{n \rightarrow + \infty}\frac{1}{|\dT_n|} \sum_{k \in\dT_{n,p}}\dX_k\dX_k^t = \ell,
\hspace{1cm}\text{a.s.}
\end{equation}
where the matrix $\ell$
is the unique solution of the equation
$$\ell=T+\frac{1}{2}(A\ell A^t+B\ell B^t)$$
\begin{equation*}
T= ( \sigma^2 + \overline{a^2})e_1e_1^t + \frac{1}{2} ( a_0( A \lambda e_1^t
+ e_1 \lambda^tA^t) + b_0( B \lambda e_1^t + e_1
\lambda^tB^t))
\end{equation*}
with $\overline{a^2}=(a_0^2+b_0^2)/2$.
\end{Lemma}

\noindent\textbf{Proof :} The proofs are given in Appendix~\ref{appendixA}.
\hspace{\stretch{1}}$\Box$\\

\begin{Remark}$\pt$\label{decmatl}
We shall see in Appendix A that
\begin{eqnarray*}
\ell=\sum_{k=0}^\infty\frac{1}{2^k} \sum_{C \in \{A;B\}^k}CTC^t
\end{eqnarray*}
where the notation $\{A;B\}^k$ means the set of all products of $A$ and $B$ with exactly $k$ terms.
For example, we have $\{A;B\}^0= \{\rI_p\}$, $\{A;B\}^1= \{A, B\}$, $\{A;B\}^2= \{A^2,AB, BA, B^2\}$ and so on. The cardinality of $\{A;B\}^k$ is obviously $2^k$.
\end{Remark}

\begin{Remark}$\pt$\label{casep1}
One can observe that in the special case $p=1$, 
\begin{eqnarray*}
\lim_{n \rightarrow + \infty}\frac{1}{|\dT_n|} \sum_{k \in\dT_{n}}X_k &=& 
\frac{\overline{a}}{1 -\overline{b}}
\hspace{1cm}\text{a.s.}\\
\lim_{n \rightarrow + \infty}\frac{1}{|\dT_n|} \sum_{k \in\dT_{n}}X_k^2 &=& 
\frac{\overline{a^2}+  \sigma^2 + 2\lambda\overline{ab}}{1-\overline{b^2}}
\hspace{1cm}\text{a.s.}
\end{eqnarray*}
where
\begin{equation*}
\overline{ab}= \frac{a_0a_1 + b_0b_1}{2},\hspace{1cm} \overline{b}=\frac{a_1 + b_1}{2}, \hspace{1cm}
\overline{b^2}=\frac{a_1^2 + b_1^2}{2}.
\end{equation*}
\end{Remark}

\section{Proof of Theorems~\ref{th Mn} and~\ref{thmaptheta}}
\label{section th1}
Theorem \ref{thmaptheta} is a consequence of Theorem~\ref{th Mn}. The first result of 
Theorem \ref{th Mn} is a strong law of large numbers for the martingale $(M_n)$. 
We already mentioned that the standard strong law is useless here.
This is due to the fact that the dimension of the random vector $\xi_n$ grows exponentially fast
as $2^n$. Consequently, we are led to propose a new strong law of large numbers for $(M_n)$,
adapted to our framework.\\

\noindent\textbf{Proof of result~(\ref{th Mn1}) of Theorem \ref{th Mn} :}
For all $n\geq p$, let $\cV_{n}=M_n^t\Sigma_{n-1}^{-1}M_n$ where we recall that 
$\Sigma_n=\rI_{2} \otimes S_{n}$,
so that $\Sigma_n^{-1}=\rI_{2} \otimes S_{n}^{-1}$.  First of all, we have
\begin{eqnarray*}
\cV_{n+1}
&\!=\!& M_{n+1}^t\Sigma_n^{-1}M_{n+1}=(M_n+\Delta M_{n+1})^t\Sigma_n^{-1}(M_n+\Delta M_{n+1}),\\
&\!=\!&M_n^t\Sigma_n^{-1}M_n+2M_n^t\Sigma_n^{-1}\Delta M_{n+1}+\Delta M_{n+1}^t\Sigma_n^{-1}\Delta M_{n+1},\\
&\!=\!&\cV_n\!-\!M_n^{t}(\Sigma_{n-1}^{-1}\!-\!\Sigma_n^{-1})M_n\!+\!
2M_n^t\Sigma_n^{-1}\Delta M_{n+1}\!+\!\Delta M_{n+1}^t\Sigma_n^{-1}\Delta M_{n+1}.
\end{eqnarray*}
By summing over this identity, we obtain the main decomposition
\begin{equation}
\label{maindecomart}
\cV_{n+1}+\cA_n=\cV_p+\cB_{n+1}+\cW_{n+1},
\end{equation}
where
\begin{eqnarray*}
\cA_n=\sum_{k=p}^n M_k^{t}(\Sigma_{k-1}^{-1}-\Sigma_k^{-1})M_k,
\end{eqnarray*}
\vspace{-0.8cm}
\begin{eqnarray*}
\cB_{n+1}=2\sum_{k=p}^n M_k^t\Sigma_k^{-1}\Delta M_{k+1}\hspace{0.6cm}\text{and}\hspace{0.6cm}
\cW_{n+1}=\sum_{k=p}^n \Delta M_{k+1}^t\Sigma_k^{-1}\Delta M_{k+1}.
\end{eqnarray*}

\noindent The asymptotic behavior of the left-hand side of (\ref{maindecomart}) is as follows.
\begin{Lemma}\label{lem lim V+A}
Assume that $(\veps_n)$ satisfies \emph{(\textbf{H.1})} to \emph{(\textbf{H.3})}. Then, we have
\begin{equation}
\label{cvgVA}
\lim_{n \rightarrow + \infty}\frac{\cV_{n+1}+\cA_n}{n} = (p+1)\sigma^2
\hspace{1cm}\text{a.s.}
\end{equation}
\end{Lemma}

\noindent\textbf{Proof :} The proof is given in Appendix B. It relies
on the Riccation equation associated to $(S_n)$ and the
strong law of large numbers for $(\cW_n)$. \hspace{\stretch{1}}$\Box$\\

\noindent Since $(\cV_{n})$ and $(\cA_n)$ are two sequences of positive real numbers, we infer from 
Lemma~\ref{lem lim V+A} that $\cV_{n+1}=\cO(n)$ a.s. which ends the proof of (\ref{th Mn1}).\hspace{\stretch{1}}$\Box$\\

\noindent\textbf{Proof of result~(\ref{aptheta1}) of Theorem \ref{thmaptheta}:}
It clearly follows from (\ref{thetadiff}) that
\begin{equation*}
\cV_n=(\wh{\theta}_n-\theta)^t\Sigma_{n-1}(\wh{\theta}_n-\theta).
\end{equation*}
Consequently, the asymptotic behavior of 
$\wh{\theta}_n-\theta$ is clearly related to the one of $\cV_n$. More precisely, 
we can deduce from convergence (\ref{cvgSn}) that
\begin{equation*}
\lim_{n\rightarrow\infty}\frac{\lambda_{\textrm{min}}(\Sigma_{n})}{|\mathbb{T}_{n}|}
=\lambda_{\textrm{min}}(\Lambda)>0\hspace{1cm} \text{a.s.}
\end{equation*}
since $L$ as well as $\Lambda=\rI_{2}\otimes L$ are definite positive matrices. Here $\lambda_{\textrm{min}}(\Lambda)$ stands for the smallest eigenvalue of the matrix $\Lambda$. 
Therefore, as
\begin{equation*}
\|\widehat{\theta}_{n}-\theta\|^{2}\leq \frac{\cV_n}{\lambda_{\textrm{min}}(\Sigma_{n-1})},
\vspace{-1ex}
\end{equation*}
we use (\ref{th Mn1}) to conclude that
\begin{equation*}
\|\widehat{\theta}_{n}-\theta\|^{2}= \cO \left(\frac{n}{|\dT_{n-1}|} \right)=
\cO \left(\frac{\log |\dT_{n-1}|}{|\dT_{n-1}|} \right)
\hspace{1cm}\text{a.s.}
\end{equation*}
which completes the proof of (\ref{aptheta1}). \hspace{\stretch{1}}$\Box$\\

We now turn to the proof of the quadratic strong law. To this end, we need a sharper estimate of the asymptotic behavior of $(\cV_n)$.

\begin{Lemma}\label{lemlimM}
Assume that $(\veps_n)$ satisfies \emph{(\textbf{H.1})} to \emph{(\textbf{H.3})}.
Then, we have for all $\delta>1/2$, 
\begin{equation}
\label{th Mn2}
\parallel M_n \parallel^2=o(|\dT_{n-1}| n^\delta)
\hspace{1cm}\text{a.s.}
\end{equation}
\end{Lemma}

\noindent\textbf{Proof :} The proof is given in Appendix C. \hspace{\stretch{1}}$\Box$\\

A direct application of Lemma~\ref{lemlimM} ensures that $\cV_n=o(n^{\delta})$ a.s. 
for all $\delta>1/2$. Hence, Lemma~\ref{lem lim V+A} immediately leads to the following result.

\begin{Corollary}\label{corlimA}
Assume that $(\veps_n)$ satisfies \emph{(\textbf{H.1})} to
\emph{(\textbf{H.3})}. Then, we have
\begin{equation}
\label{cvgAn}
\lim_{n \rightarrow + \infty}\frac{\cA_n}{n} = (p+1) \sigma^2
\hspace{1cm}\text{a.s.}
\end{equation}
\end{Corollary}

\noindent\textbf{Proof of result~(\ref{th Mn3}) of Theorem \ref{th Mn}:} First of all, $\cA_n$ may be rewritten as
\begin{equation*}
\cA_n=\sum_{k=p}^n M_k^{t}(\Sigma_{k-1}^{-1}-\Sigma_k^{-1})M_k
= \sum_{k=p}^nM_k^{t}\Sigma_{k-1}^{-1/2}\Delta_k \Sigma_{k-1}^{-1/2}M_k
\end{equation*}
where $\Delta_n=\rI_{2(p+1)} - \Sigma_{n-1}^{1/2}\Sigma_n^{-1} \Sigma_{n-1}^{1/2}$.
In addition, via Proposition~\ref{mainlemma}
\begin{equation}
\label{cvgSigman}
\lim_{n\rightarrow\infty}\frac{\Sigma_{n}}{|\mathbb{T}_{n}|}=\Lambda
\hspace{1cm}\text{a.s.}
\end{equation}
which implies that
\begin{equation}
\label{cvgDeltan}
\lim_{n\rightarrow\infty}\Delta_n=\frac{1}{2}\rI_{2(p+1)}
\hspace{1cm}\text{a.s.}
\end{equation}
Furthermore, it follows from Corollary~\ref{corlimA} that $\cA_{n}=\cO(n)$ a.s. 
Hence, we deduce from (\ref{cvgSigman}) and (\ref{cvgDeltan})
that
\begin{equation}\label{ApproxAn}
\frac{\cA_n}{n}= \left( \frac{1}{2n}\sum_{k=p}^nM_k^{t}\Sigma_{k-1}^{-1} M_k \right)+ o(1)
\hspace{1cm}\text{a.s.} \\
\end{equation}
and convergence~(\ref{th Mn3}) directly follows from Corollary~\ref{corlimA}.\hspace{\stretch{1}}$
\Box$\\

We are now in position to prove the QSL.\\

\noindent\textbf{Proof of result~(\ref{aptheta2}) of Theorem \ref{thmaptheta}:}
The QSL is a direct consequence of (\ref{th Mn3}) together with the fact that
$\wh{\theta}_n-\theta=\Sigma_{n-1}^{-1}M_n$. Indeed, we have
\begin{eqnarray*}
\frac{1}{n}\sum_{k=p}^nM_k^{t}\Sigma_{k-1}^{-1} M_k 
&=&\frac{1}{n}\sum_{k=p}^n (\wh{\theta}_k-\theta)^{t}\Sigma_{k-1}(\wh{\theta}_k-\theta)
 \\
&=&\frac{1}{n}\sum_{k=p}^n |\dT_{k-1}|(\wh{\theta}_k-\theta)^{t}
\frac{\Sigma_{k-1}}{|\dT_{k-1}|}(\wh{\theta}_k-\theta)
 \\
&=& \frac{1}{n}\sum_{k=p}^n |\dT_{k-1}|(\wh{\theta}_k-\theta)^{t}
\Lambda(\wh{\theta}_k-\theta) + o(1)
\hspace{1cm}\text{a.s.}
\end{eqnarray*}
which completes the proof of Theorem \ref{thmaptheta}.\hspace{\stretch{1}}$
\Box$\\

\section{Proof of Theorem~\ref{thmapsigmarho}}
\label{section th2}
The almost sure convergence of $\wh{\sigma}^2_n$ and $\wh{\rho}_n$ is strongly related
to that of $\wh{V}_n-V_n$.\\

\noindent\textbf{Proof of result~(\ref{apsigma1}) of Theorem~\ref{thmapsigmarho}:} We need to prove that
\begin{equation}
\label{cvgVchap}
\lim_{n\rightarrow\infty}\frac{1}{n}\sum_{k\in\mathbb{T}_{n-1,p-1}}\|\wh{V}_k-V_k\|^2=2(p+1)\sigma^2\hspace{1cm}\text{a.s.}
\end{equation}
Once again, we are searching for a link between the sum of $\|\wh{V}_n-V_n\|$ and the processes $(\cA_n)$ and $(\cV_n)$ whose convergence properties were previously investigated. For all $n\geq p$, we have
\begin{eqnarray*}
\sum_{k\in\mathbb{G}_n}\|\wh{V}_k-V_k\|^2&=&\sum_{k\in\mathbb{G}_n}(\wh{\veps}_{2k}-\veps_{2k})^2+(\wh{\veps}_{2k+1}-\veps_{2k+1})^2,\\
&=&(\wh{\theta}_n-\theta)^t\Psi_n\Psi_n^t(\wh{\theta}_n-\theta),\\
&=&M_n^t\Sigma_{n-1}^{-1}\Psi_n\Psi_n^t\Sigma_{n-1}^{-1}M_n,\\
&=&M_n^t\Sigma_{n-1}^{-1/2}\Delta_n\Sigma_{n-1}^{-1/2}M_n,
\end{eqnarray*}
where\\
$$
\vspace{2ex}
\Delta_n=\Sigma_{n-1}^{-1/2}\Psi_n\Psi_n^t\Sigma_{n-1}^{-1/2}=
\Sigma_{n-1}^{-1/2}(\Sigma_n-\Sigma_{n-1})\Sigma_{n-1}^{-1/2}.
$$
Now, we can deduce from convergence (\ref{cvgSigman}) that
\begin{equation*}
\lim_{n\rightarrow\infty}\Delta_n=\rI_{2(p+1)}\hspace{1cm}\text{a.s.}
\end{equation*}
which implies that
\begin{equation*}
\sum_{k\in\mathbb{G}_n}\|\wh{V}_k-V_k\|^2=M_n^t\Sigma_{n-1}^{-1}M_n\Bigl(1+o(1)\Bigr)\hspace{1cm}\text{a.s.}
\end{equation*}
Therefore, we can conclude via convergence (\ref{th Mn3}) that
\begin{equation*}
\lim_{n\rightarrow\infty}\frac{1}{n}\sum_{k\in\mathbb{T}_{n-1,p-1}}\|\wh{V}_k-V_k\|^2=
\lim_{n\rightarrow\infty}\frac{1}{n}\sum_{k=p}^n M_k^{t}\Sigma_{k-1}^{-1}M_k=2(p+1)\sigma^2\hspace{1cm}\text{a.s.}
\end{equation*}

\noindent\textbf{Proof of result~(\ref{apsigma2}) of Theorem~\ref{thmapsigmarho}:} First of all,
\begin{eqnarray*}
\wh{\sigma}^2_n-{\sigma}^2_n&=&\frac{1}{2|\mathbb{T}_{n-1}|}
\sum_{k\in\mathbb{T}_{n-1,p-1}}\big(\|\wh{V}_k\|^2-\|{V}_k\|^2\big),\\
&=&\frac{1}{2|\mathbb{T}_{n-1}|}\sum_{k\in\mathbb{T}_{n-1,p-1}}\big(\|\wh{V}_k-{V}_k\|^2+2(\wh{V}_k-{V}_k)^tV_k\big).
\end{eqnarray*}
Set 
\begin{equation*}
P_n=\sum_{k\in\mathbb{T}_{n-1,p-1}}(\wh{V}_k-{V}_k)^tV_k=
\sum_{k=p}^{n}\sum_{i\in\mathbb{G}_{k-1}}(\wh{V}_i-{V}_i)^tV_i.
\end{equation*}
We clearly have
\begin{equation*}
\Delta P_{n+1}=P_{n+1}-P_n=\sum_{k \in \dG_{n}}(\wh{V}_k-{V}_k)^tV_k.
\end{equation*}
One can observe that for all $k \in \mathbb{G}_{n}$, $\wh{V}_k-{V}_k=(\rI_2\otimes Y_k)^t(\theta-\wh{\theta}_n)$ which implies that $\wh{V}_k-{V}_k$ is $\cF_n$-measurable. 
Consequently, $(P_n)$ is a real martingale transform.
Hence, we can deduce from the strong law of large numbers for martingale
transforms given in Theorem 1.3.24 of \cite{Duflo97} together with (\ref{cvgVchap}) that 
\begin{equation*}
P_n=o\left(\sum_{k\in\dT_{n-1,p-1}}||\wh{V}_k-{V}_k)||^2\right)=o(n)\hspace{1cm}\text{a.s.}
\end{equation*}
It ensures once again via convergence (\ref{cvgVchap}) that
\begin{equation*}
\lim_{n\rightarrow\infty}\frac{|\dT_{n}|}{n}
(\wh{\sigma}^2_n-{\sigma}^2_n)
=\lim_{n\rightarrow\infty}\frac{1}{n}\sum_{k\in\mathbb{T}_{n-1,p-1}}\|\wh{V}_k-V_k\|^2=
2(p+1)\sigma^2\hspace{1cm}\text{a.s.}
\end{equation*}
We now turn to the study of the covariance estimator $\wh{\rho}_n$. We have
\begin{eqnarray*}
\wh{\rho}_n-{\rho}_n&=&\frac{1}{|\dT_{n-1}|}\sum_{k\in\mathbb{T}_{n-1,p-1}}(\wh{\veps}_{2k}\wh{\veps}_{2k+1}-\veps_{2k}\veps_{2k+1}),\\
&=&\frac{1}{|\dT_{n-1}|}\sum_{k\in\mathbb{T}_{n-1,p-1}}(\wh{\veps}_{2k}-\veps_{2k})(\wh{\veps}_{2k+1}-\veps_{2k+1})
+\frac{1}{|\dT_{n-1}|}Q_n,
\end{eqnarray*}
where
$$Q_n=\!\!\!\!\sum_{k\in\mathbb{T}_{n-1,p-1}}\!\!\!\!
(\wh{\veps}_{2k}-\veps_{2k})\veps_{2k+1}+(\wh{\veps}_{2k+1}-\veps_{2k+1})\veps_{2k}
=\!\!\!\!\sum_{k\in\mathbb{T}_{n-1,p-1}}\!\!\!\!(\wh{V}_k-{V}_k)^t\rJ_2V_k$$
with
\begin{equation*}
\rJ_2=\left(\begin{array}{cc}0&1\\1&0\end{array}\right).
\end{equation*}
Moreover, one can observe that $\rJ_2\Gamma \rJ_2= \Gamma$. Hence, as before,
$(Q_n)$ is a real martingale transform satisfying
\begin{equation*}
Q_n=o\left(\sum_{k\in\dT_{n-1,p-1}}||\wh{V}_k-{V}_k)||^2\right)=o(n)\hspace{1cm}\text{a.s.}
\end{equation*}
We will see in Appendix~\ref{appendixD} that
\begin{equation}
\label{cvgVJchap}
\lim_{n\rightarrow\infty}\frac{1}{n}
\sum_{k\in\dT_{n-1,p-1}}(\wh{\varepsilon}_{2k}-\varepsilon_{2k})(\wh{\varepsilon}_{2k+1}-\varepsilon_{2k+1})
=(p+1)\rho
\hspace{1cm}\text{a.s.}
\end{equation}
Finally, we find from (\ref{cvgVJchap}) that
\begin{equation*}
\lim_{n\rightarrow\infty}\frac{|\dT_{n}|}{n}
(\wh{\rho}_n-{\rho}_n)=
2(p+1)\rho\hspace{1cm}\text{a.s.}
\end{equation*}
which completes the proof of Theorem~\ref{thmapsigmarho}.\hspace{\stretch{1}}$
\Box$\\

\section{Proof of Theorem~\ref{thmCLT}}
\label{section th3}
In order to prove the CLT for the BAR($p$) estimators, 
we will use the central limit theorem for martingale difference
sequences given in Propositions~7.8 and 7.9 of Hamilton \cite{Ham94}.

\begin{Proposition}\label{prop Hamilton}
Assume that $(W_n)$ is a vector martingale difference sequence satisfying
\begin{itemize}
\item[(a)] $\hspace{-1ex}$ For all $n\!\geq\! 1$, $\dE[W_nW_n^t]\!=\!\Omega_n$ where $\Omega_n\!$ is
a positive definite matrix and
$$
\lim_{n\rightarrow\infty}\frac{1}{n} \sum_{k=1}^n\Omega_k=\Omega
$$
where $\Omega$ is also a positive definite matrix.
\item[(b)] For all $n\geq 1$ and for all $i,j,k,l$, $\dE[W_{in}W_{jn}W_{kn}W_{ln}]<\infty$ where 
$W_{in}$ is the $i$th element of the vector $W_n$.
\item[(c)] 
$$\frac{1}{n}\sum_{k=1}^n W_kW_k^t\build{\longrightarrow}_{}^{{\mbox{\calcal P}}}\Omega.
$$
\end{itemize}
Then, we have the central limit theorem
\begin{equation*}
\frac{1}{\sqrt{n}}\sum_{k=1}^n W_k \build{\longrightarrow}_{}^{{\mbox{\calcal L}}}
\cN(0,\Omega).
\end{equation*}
\end{Proposition}

We wish to point out that for BAR$(p)$ processes, it seems impossible to make use of
the standard CLT for martingales. This is due to the fact that Lindeberg's condition
is not satisfied in our framework. Moreover, as the size of $(\xi_n)$ doubles at each generation, 
it is also impossible to check condition $(c)$. To overcome this problem, we simply change the filtration.
Instead of using the generation-wise filtration, we will use the sister pair-wise one. Let
\begin{equation*}
\mathcal{G}_n=\sigma\{X_1,\ (X_{2k}, X_{2k+1}),\ 1\leq k\leq n\}
\end{equation*}
be the $\sigma$-algebra generated by all pairs of individuals up to the offspring of individual $n$. Hence $(\veps_{2n}, \veps_{2n+1})$ is $\mathcal{G}_n$-measurable. Note that $\mathcal{G}_n$ is also the $\sigma$-algebra generated by, on the one hand, all the past generations up to that of individual $n$, i.e. the $r_n$th generation, and, on the other hand, all pairs of the $(r_n+1)$th generation with ancestors less than or equal to $n$. In short,
\begin{equation*}
\mathcal{G}_n=\sigma\Bigl(\mathcal{F}_{r_n}\cup\{(X_{2k}, X_{2k+1}), \ k\in\mathbb{G}_{r_n}, \ k\leq n\}\Bigr).
\end{equation*}
Therefore, {(\textbf{H.2})} implies that the processes $(\veps_{2n}, \dX_n\veps_{2n}, \veps_{2n+1}, \dX_n\veps_{2n+1})^t$, $(\veps_{2n}^2+\veps_{2n+1}^2-2\sigma^2)$ and $(\veps_{2n}\veps_{2n+1}-\rho)$ are $\mathcal{G}_n$-martingales.\\

\noindent\textbf{Proof of result~(\ref{th Mn4}) of Theorem~\ref{th Mn}:} First, recall that
$Y_n=(1,\dX_n)^t$. We apply Propositions~\ref{prop Hamilton} 
to the $\mathcal{G}_n$-martingale difference sequence
$(D_n)$ given by
\begin{equation*}
D_n=\mathrm{vec}(Y_nV_n^t)=\left( \begin{array}{cccc}
\veps_{2n}  \\
\dX_n\veps_{2n} \\
\veps_{2n+1} \\
\dX_n\veps_{2n+1}
\end{array}\right).
\end{equation*}
We clearly have
\begin{equation*}
D_nD_n^t=\left(
\begin{array}{ll}
\varepsilon_{2n}^2 & \varepsilon_{2n}\varepsilon_{2n+1}\\
\varepsilon_{2n+1}\varepsilon_{2n}&\varepsilon_{2n+1}^2
\end{array}
\right)\otimes Y_nY_n^t.
\end{equation*}
Hence, it follows from (\textbf{H.1}) and (\textbf{H.2}) that 
\begin{equation*}
\mathbb{E}[D_nD_n^t]=\Gamma\otimes\dE[Y_nY_n^t].
\end{equation*}
Moreover, we can show by a slight change in the proof of 
Lemmas~\ref{lemlimsumX} and \ref{lemlimsumXX} that
\begin{equation*}
\lim_{n\rightarrow\infty}\frac{1}{|\mathbb{T}_{n}|}
\sum_{k\in\mathbb{T}_{n-1,p-1}}\mathbb{E}[D_kD_k^t]=
\Gamma\otimes \lim_{n\rightarrow\infty}
\frac{1}{|\mathbb{T}_{n}|}\mathbb{E}[S_{n}]=\Gamma\otimes L,
\end{equation*}
which is positive definite, so that condition $(a)$ holds. Condition $(b)$ also clearly holds under {(\textbf{H.3})}. We now turn to condition $(c)$.
We have
\begin{equation*}
\sum_{k\in \mathbb{T}_{n-1,p-1}} D_kD_k^t=\Gamma\otimes S_{n}+R_n
\end{equation*}
where
\begin{equation*}
R_n=\sum_{k\in \mathbb{T}_{n-1,p-1}}\left(
\begin{array}{ll}
\varepsilon_{2k}^2-\sigma^2 & \varepsilon_{2k}\varepsilon_{2k+1}-\rho\\
\varepsilon_{2k+1}\varepsilon_{2k}-\rho&\varepsilon_{2k+1}^2-\sigma^2
\end{array}
\right)\otimes Y_k Y_k^t.
\end{equation*}
Under (\textbf{H.1}) to (\textbf{H.5}), we can show that $(R_n)$ is a 
martingale transform. Moreover, we can prove that $R_n=o(n)$ a.s. using Lemma \ref{lemsumX4}
and similar calculations as in Appendix B where a more complicated martingale
transform $(K_n)$ is studied. Consequently, condition $(c)$ also holds and we can conclude that
\begin{equation}
\label{CLTDn}
\frac{1}{\sqrt{|\dT_{n-1}|}}\sum_{k \in \mathbb{T}_{n-1,p-1}}D_k
=\frac{1}{\sqrt{|\dT_{n-1}|}}M_n
\build{\longrightarrow}_{}^{{\mbox{\calcal L}}}
\cN(0,\Gamma\otimes L).
\end{equation}

\noindent\textbf{Proof of result~(\ref{CLTtheta}) of Theorem~\ref{thmCLT}:} 
We deduce from (\ref{thetadiff}) that
\begin{equation*}
{\sqrt{|\dT_{n-1}|}}(\widehat{\theta}_n-\theta)=|\dT_{n-1}|\Sigma_{n-1}^{-1}
\frac{M_n}{\sqrt{|\dT_{n-1}|}}.
\end{equation*}
Hence, (\ref{CLTtheta}) directly follows from 
(\ref{th Mn4}) and convergence (\ref{cvgSigman}) together with Slutsky's Lemma.
\hspace{\stretch{1}}$\Box$\\

\noindent\textbf{Proof of results~(\ref{CLTsigma}) and (\ref{CLTrho}) of Theorem~\ref{thmCLT}:} On the one hand, we apply Propositions~\ref{prop Hamilton} to the $\mathcal{G}_n$-martingale difference sequence $(v_n)$ defined by
\begin{equation*}
v_n=\veps_{2n}^2+\veps_{2n+1}^2-2\sigma^2.
\end{equation*}
Under (\textbf{H.4}), one has
$\mathbb{E}[v_n^2]=2\tau^4-4\sigma^4+2\nu^2$
which ensures that
\begin{equation*}
\lim_{n\rightarrow\infty} \frac{1}{|\mathbb{T}_{n}|}\sum_{k \in \mathbb{T}_{n,p-1}} \mathbb{E}[v_k^2]=2\tau^4-4\sigma^4+2\nu^2>0.
\end{equation*}
Hence, condition $(a)$ holds. Once again, condition $(b)$ clearly holds under (\textbf{H.5}), and Lemma~\ref{lemLLEps4TCL} together with Remark~\ref{rqCVeps4} imply condition $(c)$,
\begin{equation*}
\lim_{n\rightarrow\infty}\frac{1}{|\mathbb{T}_{n}|}\sum_{k \in \mathbb{T}_{n,p-1}} v_k^2=2\tau^4-4\sigma^4+2\nu^2\hspace{1cm}\text{a.s.}
\end{equation*}
Therefore, we obtain that
\begin{equation}
\label{CLTsigman1}
\frac{1}{\sqrt{|\dT_{n-1}|}}\!\sum_{k \in \mathbb{T}_{n-1,p-1}}\!\!v_k=2\sqrt{|\dT_{n-1}|}(\sigma_n^2-\sigma^2)
\build{\longrightarrow}_{}^{{\mbox{\calcal L}}}
\cN(0,2\tau^4-4\sigma^4+2\nu^2).
\end{equation}
Furthermore, we infer from (\ref{apsigma2}) that
\begin{equation}
\label{CLTsigman2}
\lim_{n\rightarrow\infty}\sqrt{|\dT_{n-1}|}(\wh{\sigma}_n^2-\sigma_n^2)=0\hspace{1cm}\text{a.s.}
\end{equation}
Finally, (\ref{CLTsigman1}) and (\ref{CLTsigman2}) imply (\ref{CLTsigma}). On the other hand, we
apply again Proposition~\ref{prop Hamilton} to the $\mathcal{G}_n$-martingale difference sequence
$(w_n)$ given by
\begin{equation*}
w_n=\veps_{2n}\veps_{2n+1}-\rho.
\end{equation*}
Under (\textbf{H.4}), one has $\mathbb{E}[w_n^2]=\nu^2-\rho^2$
which implies that condition $(a)$ holds since
\begin{equation*}
\lim_{n\rightarrow\infty} \frac{1}{|\mathbb{T}_{n}|}\sum_{k \in\mathbb{T}_{n,p-1}} \mathbb{E}[w_k^2]=\nu^2-\rho^2>0.
\end{equation*}
Once again, condition $(b)$ clearly holds under (\textbf{H.5}), and Lemmas~\ref{lemLLNeps}
and \ref{lemLLEps4TCL} yield condition $(c)$,
\begin{equation*}
\lim_{n\rightarrow\infty}\frac{1}{|\mathbb{T}_{n}|}\sum_{k \in\mathbb{T}_{n,p-1}} w_k^2=\nu^2-\rho^2\hspace{1cm}\text{a.s.}
\end{equation*}
Consequently, we obtain that
\begin{equation}
\label{CLTrhon1}
\frac{1}{\sqrt{|\dT_{n-1}|}}\sum_{k \in \mathbb{T}_{n-1,p-1}}w_k=\sqrt{|\dT_{n-1}|}(\rho_n-\rho)
\build{\longrightarrow}_{}^{{\mbox{\calcal L}}}
\cN(0,\nu^2-\rho^2).
\end{equation}
Furthermore, we infer from (\ref{aprho2}) that
\begin{equation}
\label{CLTrhon2}
\lim_{n\rightarrow\infty}\sqrt{|\dT_{n-1}|}(\wh{\rho}_n-\rho_n)=0\hspace{1cm}\text{a.s.}
\end{equation}
Finally, (\ref{CLTrho}) follows from (\ref{CLTrhon1}) and (\ref{CLTrhon2}) which completes
the proof of Theorem~\ref{thmCLT}.
\hspace{\stretch{1}}$\Box$\\

\appendix
\section{}
\begin{center}
{\bf Laws of large numbers for the BAR process}
\end{center}
\label{appendixA}

We start with some technical Lemmas we make repeatedly use of, the well-known Kronecker's Lemma given in Lemma 1.3.14 of
\cite{Duflo97} together with some related results.

\begin{Lemma}$\pt$
\label{lemKro}
Let $(\alpha_n)$ be a sequence of positive real numbers increasing to infinity. In addition, let $(x_n)$ be a sequence of real numbers such that
$$\sum_{n=0}^{\infty}\frac{|x_n|}{\alpha_n}<+\infty.$$
Then, one has
\begin{equation*}
\lim_{n \rightarrow \infty}\frac{1}{\alpha_n}\sum_{k=0}^{n}{x_k}=0.
\end{equation*}
\end{Lemma}

\begin{Lemma}$\pt$
\label{lem tree-generation}
Let $(x_n)$ be a sequence of real numbers. Then,
\begin{equation}
\label{equi tree-generation} 
\lim_{n\rightarrow \infty}\frac{1}{|\dT_n|}\sum_{k\in\dT_n}x_k=x
\Longleftrightarrow
\lim_{n\rightarrow \infty}\frac{1}{|\dG_n|}\sum_{k\in\dG_n}x_k=x.
\end{equation}
\end{Lemma}
 
\noindent \textbf{Proof:} First of all, recall that $|\dT_n|=2^{n+1}-1$ and $|\dG_n|=2^n$. Assume that
$$
\lim_{n\rightarrow \infty}\frac{1}{|\dT_n|}\sum_{k\in\dT_n}x_k=x.
$$
We have the decomposition,
$$ 
\sum_{k\in\dT_n}x_k = \sum_{k\in\dT_{n-1}}x_k + \sum_{k\in\dG_n}x_k.
$$
Consequently,
\begin{eqnarray*}
\lim_{n\rightarrow \infty}
\frac{1}{|\dG_n|}\sum_{k\in\dG_n}x_k&=&
\lim_{n\rightarrow \infty}\frac{2}
{|\dT_n|+1}\sum_{k\in\dT_n}x_k-
\lim_{n\rightarrow \infty}
\frac{1}{|\dT_{n-1}|+1}\sum_{k\in\dT_{n-1}}x_k,\\
&=&2x-x\ =\ x.
\end{eqnarray*}
Conversely, suppose that
$$
\lim_{n\rightarrow \infty}\frac{1}{|\dG_n|}\sum_{k\in\dG_n}x_k=x.
$$
A direct application of Toeplitz Lemma given in Lemma 2.2.13 of \cite{Duflo97})
yields
\begin{eqnarray*}
\lim_{n\rightarrow \infty}
\frac{1}{|\dT_n|}\sum_{k\in\dT_n}x_k
&=&
\lim_{n\rightarrow \infty}
\frac{1}{|\dT_n|} \sum_{k=0}^n \sum_{i\in\dG_k}x_i,\\
&=&
\lim_{n\rightarrow \infty}
\frac{1}{|\dT_n|} \sum_{k=0}^n 2^k\frac{1}{|\dG_k|}\sum_{i\in\dG_k}x_i
=x.
\end{eqnarray*}

\begin{Lemma}$\pt$
\label{lemtoepbar}
Let $(A_n)$ be a sequence of real-valued matrices such that 
$\sum_{n=0}^{\infty}\|A_n\|<\infty$ and
$$\lim_{n \rightarrow \infty} \sum_{k=0}^nA_k=A.$$
In addition, let $(X_n)$ be a sequence of real-valued vectors which converges to a limiting value $X$. 
Then,
\begin{equation}
\label{cvgtoepbar}
\lim_{n \rightarrow \infty}{\sum_{k=0}^n A_{n-k}X_k}=AX.
\end{equation}
\end{Lemma}

\noindent\textbf{Proof:} For all $n\geq 0$, let 
$$U_n=\sum_{k=0}^n A_{n-k}X_k.$$ 
We clearly have for all integer $n_0$ with $1\leq n_0 <n$,
\begin{eqnarray*}
\|U_n-AX\|
&=&\Big\|\sum_{k=0}^nA_{n-k}X_k-\sum_{k=0}^nA_kX-\sum_{k=n+1}^{\infty}A_kX\Big\|,\\
&\leq&\sum_{k=0}^n\|A_{n-k}\|\|X_k-X\|+\sum_{k=n+1}^{\infty}\|A_k\|\|X\|,
\end{eqnarray*}
\vspace{-5ex}
\begin{eqnarray*}
\leq\sum_{k=0}^{n_0}\|A_{n-k}\|\|X_k-X\|+\!\!\sum_{k=n_0+1}^n\|A_{n-k}\|\|X_k-X\|+\!\!\sum_{k=n+1}^{\infty}\|A_k\|\|X\|.
\end{eqnarray*}
We assume that $(X_n)$ converges to a limiting value $X$.
Consequently, we can choose $n_0$ such that for all $k>n_0$, 
$\|X_k-X\|<\varepsilon$. Moreover, one can find $M>0$ such that for all 
$k \geq 0$, $\|X_k-X\|\leq M$ and $\|X\|\leq M$. Therefore, we
obtain that
\begin{equation*}
\|U_n-AX\|\leq (n_0+1) M \!\!
\sup_{k \geq n-n_0}\!\!\|A_k\|+\varepsilon \sum_{k=n_0+1}^n\|A_{n-k}\|+M\!\sum_{k=n+1}^{\infty}\|A_k\|.
\end{equation*}
On the one hand 
$$\sup_{k\geq n-n_0}\|A_k\| 
\hspace{1cm}\text{and} \hspace{1cm}
\sum_{k=n+1}^{\infty}\|A_k\|$$ 
both converge to $0$ as $n$ tends to infinity. On the other hand, 
$$\sum_{k=n_0+1}^n\|A_{n-k}\|\leq \sum_{n=0}^{\infty}\|A_n\|<\infty.$$ 
Consequently, $\|U_n-AX\|$ goes to $0$ as $n$ goes to infinity, as expected.
\hspace{\stretch{1}} $\Box$\\

\begin{Lemma}$\pt$
\label{lempartX2} Let $(T_n)$ be a convergent sequence of
real-valued matrices with limiting value $T$. Then,
$$\lim_{n \rightarrow \infty} \sum_{k=0}^n\frac{1}{2^k} \sum_{C \in \{A;B\}^k}CT_{n-k}C^t =
\ell $$ 
where the matrix 
$$ \ell=\sum_{k=0}^\infty\frac{1}{2^k} \sum_{C \in \{A;B\}^k}CTC^t$$ 
is the unique solution of the equation
\begin{equation}
\label{eqrecl}
\ell = T + \frac{1}{2}(A \ell A^t + B \ell B^t).
\end{equation}
\end{Lemma}

\noindent\textbf{Proof: } First of all, recall that
$\beta=\max\{\|A\|, \|B\|\}<1$. The cardinality of $\{A;B\}^k$ is obviously $2^k$.
Consequently, if
$$
U_n= \sum_{k=0}^n\frac{1}{2^k}\sum_{C \in \{A;B\}^k}C(T_{n-k}-T)C^t,
$$
it is not hard to see that
$$
\|U_n \| \leq  \sum_{k=0}^n\frac{1}{2^k}\times 2^k \beta^{2k} \Big\|T_{n-k}-T\Big\|
=\sum_{k=0}^n \beta^{2(n-k)} \Big\|T_{k}-T\Big\|.
$$
Hence, $(U_n)$ converges to zero which completes the proof of Lemma~\ref{lempartX2}.
\hspace{\stretch{1}} $\Box$\\

We now return to the BAR process. We first need an estimate of the sum of the $\|\dX_n\|^2$ before being able to investigate the limits.
\begin{Lemma}$\pt$
\label{lemSumX}
Assume that $(\veps_n)$ satisfies \emph{(\textbf{H.1})} to \emph{(\textbf{H.3})}. Then, we have
\begin{equation}
\label{sumX2}
\sum_{k \in \dT_{n,p}}\|\dX_k\|^2 =\cO (|\dT_n|)
\hspace{1cm}\text{a.s.}
\end{equation}
\end{Lemma}

\noindent\textbf{Proof :} In all the sequel, for all $n\geq 2^{p-1}$, denote $A_{2n} = A$ and $A_{2n+1} = B$. It follows from a recursive application of relation~(\ref{defbarmatrix}) that
for all $n\geq 2^{p-1}$
\begin{equation}\label{eq X itere}
\dX_n = \Big(\prod_{k=0}^{r_n-p}A_{[\frac{n}{2^k}]}\Big)\dX_{[\frac{n}{2^{r_n-p+1}}]}+ \sum_{k=0}^{r_n-p}
\Big(\prod_{i=0}^{k-1}A_{[\frac{n}{2^i}]}\Big)\eta_{[\frac{n}{2^k}] }
\end{equation}
with the convention that an empty product equals $1$. Then, 
we can deduce from Cauchy-Schwarz inequality that for all $n\geq 2^{p-1}$
\begin{eqnarray*}
\Big\|\dX_n-\big(\prod_{k=0}^{r_n-p}A_{[\frac{n}{2^k}]}\big)\dX_{[\frac{n}{2^{r_n-p+1}}]} \Big\|^2
&=& \left\| \sum_{k=0}^{r_n-p}
\Big(\prod_{i=0}^{k-1}A_{[\frac{n}{2^i}]}\Big)\eta_{[\frac{n}{2^k}] }\right\|^2\\
&\leq&\left( \sum_{k=0}^{r_n-p}
\Big(\prod_{i=0}^{k-1}\|A_{[\frac{n}{2^i}]}\|\Big)\|\eta_{[\frac{n}{2^k}]}\|\right)^2\\
& \leq & \left( \sum_{k=0}^{r_n-p}\beta^k\big\|\eta_{[\frac{n}{2^k}] }\big\|\right)^2\\
&\leq & \left( \sum_{k=0}^{r_n-p}\beta^k\right)\left( \sum_{k=0}^{r_n-p}\beta^k\|\eta_{[\frac{n}{2^k}] }\|^2\right)\\
&\leq & \frac{1}{1-\beta}\left( \sum_{k=0}^{r_n-p}\beta^k\|\eta_{[\frac{n}{2^k}] }\|^2\right).
\end{eqnarray*}
Hence, we obtain that for all $n\geq 2^p$,
\begin {eqnarray*}
\|\dX_n\|^2 & = & \left\|\dX_n-\Big(\prod_{k=0}^{r_n-p}A_{[\frac{n}{2^k}]}\Big)\dX_{[\frac{n}{2^{r_n-p+1}}]} + \Big(\prod_{k=0}^{r_n-p}A_{[\frac{n}{2^k}]}\Big)\dX_{[\frac{n}{2^{r_n-p+1}}]}\right\|^2 \\
& \leq & \frac{2}{1-\beta}\left( \sum_{k=0}^{r_n-p}\beta^k\|\eta_{[\frac{n}{2^k}]}\|^2\right) + 2\beta^{2(r_n-p+1)}\|\dX_{[\frac{n}{2^{r_n-p+1}}]}\|^2.
\end{eqnarray*}
Denote $\alpha=\max\{|a_0|,|b_0|\}$ and $\overline{X}_1=\max\{\|X_k\|, k\leq2^{p-1}\}$. Summing up over the sub-tree $\dT_{n,p}$, we find that
\begin{eqnarray}
\sum_{k \in \dT_{n,p}}\!\! \|\dX_k\|^2\!\!
& \!\leq \! & \!\!\sum_{k \in \dT_{n,p}}\!\!
\frac{2}{1-\beta}\left( \sum_{i=0}^{r_k-p}\beta^i\|\eta_{[\frac{k}{2^i}]}\|^2\right)
+\!\! \sum_{k \in \dT_{n,p}}\!\!2\beta^{2(r_k-p+1)}\|\dX_{[\frac{k}{2^{r_k-p+1}}]}\|^2 \nonumber \\
& \!\leq\! & \!\!\frac{4}{1-\beta}\!\!\sum_{k \in \dT_{n,p}}\!\!
\sum_{i=0}^{r_k-p}\beta^i(\alpha^2 + \veps^2_{[\frac{k}{2^i}] })
+ \!\!\sum_{k \in \dT_{n,p}}\!\!2\beta^{2(r_k-p+1)}\|X_{[\frac{k}{2^{r_k-p+1}}]}\|^2 \nonumber  \\
& \!\leq\! & \!\!\frac{4}{1-\beta}\sum_{k \in \dT_{n,p}}
\sum_{i=0}^{r_k-p}\beta^i\veps^2_{[\frac{k}{2^i}] }
+ \frac{4\alpha^2}{1-\beta}\sum_{k \in \dT_{n,p}}\sum_{i=0}^{r_k-p}\beta^i\nonumber\\
&&\!\!+ 2\overline{X_1}^2\sum_{k \in \dT_{n,p}}\beta^{2(r_k-p+1)}, \nonumber \\
& \!\leq\! &\!\! \frac{4P_n}{1-\beta}+\frac{4\alpha^2Q_n}{1-\beta}+2\overline{X}_1^2R_n,
\label{decosumX}
\end{eqnarray}
where
\begin{equation*}
P_n=\sum_{k \in \dT_{n,p}} \sum_{i=0}^{r_k-p}\beta^i\veps^2_{[\frac{k}{2^i}]},\ \
Q_n=\sum_{k \in \dT_{n,p}}\sum_{i=0}^{r_k-p}\beta^i,\ \
R_n=\sum_{k \in \dT_{n,p}}\!\!\beta^{2(r_k-p+1)}.
\end{equation*}
The last two terms of (\ref{decosumX}) are readily evaluated by splitting the sums generation-wise.
As a matter of fact,
\begin{equation}
\label{resBn}
Q_n=\sum_{k=p}^n\sum_{i \in \dG_k}\frac{1-\beta^k}{1-\beta} \leq \frac{1}{(1-\beta)}\sum_{k=p}^n2^k=\cO(|\dT_n|),
\end{equation}
and
\begin{equation}
\label{resCn}
R_n=\sum_{k=p}^n\sum_{i \in \dG_k}\beta^{k-p+1}\leq\sum_{k=p}^n(2\beta)^k=\cO(|\dT_n|).
\end{equation}
It remains to control the first term $P_n$. One can observe that $\veps_k$ appears in $P_n$
as many times as it has descendants up to the $n$th generation, and
its multiplicative factor for its $i$th generation descendant is $(2\beta)^i$. Hence, one has
\begin{equation*}
P_n = \sum_{k \in \dT_{n,p}} \sum_{i=0}^{n-r_k}(2\beta)^i \veps_k^2.
\end{equation*}
The evaluation of $P_n$ depends on the value of $0<\beta<1$. 
On the one hand, if $\beta= 1/2$, $P_n$ reduces to
\begin{equation*}
P_n = \sum_{k \in\dT_{n,p}}(n+1-r_k)\veps_k^2=\sum_{k=p}^n(n+1-k)\sum_{i \in \dG_k}\veps_i^2.
\end{equation*}
Hence,
\begin{equation*}
\frac{P_n}{|\dT_{n}|+1} = \sum_{k=p}^n\left(\frac{(n+1-k)}{2^{n+1-k}}\right)
\left(\frac{1}{|\dG_k|}\sum_{i \in \dG_k}\veps_i^2\right).
\end{equation*}
However, it follows from Remark~\ref{rqCVeps2k} that
\begin{equation*}
\lim_{n \rightarrow + \infty}\frac{1}{|\dG_n|}\sum_{k \in\dG_n}\veps_k^2=\sigma^2
\hspace{1cm}\text{a.s.}
\end{equation*}
In addition, we also have
$$\lim_{n \rightarrow \infty} \sum_{k=1}^n\frac{k}{2^k}=2.$$
Consequently, we infer from Lemma \ref{lemtoepbar} that
\begin{equation}
\label{cvgAn1}
\lim_{n \rightarrow + \infty}\frac{P_n}{|\dT_{n}|}=2\sigma^2
\hspace{1cm}\text{a.s.}
\end{equation}
On the other hand, if $\beta \neq 1/2$, we have
\begin{equation*}
P_n = \sum_{k \in\dT_{n,p}}\frac{1-(2\beta)^{n-r_k+1}}{1-2\beta}\veps_k^2
=\frac{1}{1-2\beta}\sum_{k=p}^n(1-(2\beta)^{n-k+1})\sum_{i \in \dG_k}\veps_i^2.
\end{equation*}
Thus,
\begin{equation*}
\frac{P_n}{|\dT_{n}|+1}=\frac{1}{1-2\beta}\sum_{k=p}^n
\left(\Big(\frac{1}{2}\Big)^{n-k+1}-\beta^{n-k+1}\right)\left(\frac{1}{|\dG_k|}\sum_{i \in \dG_k}\veps_i^2\right).
\end{equation*}
Furthermore,
$$\lim_{n \rightarrow \infty} \frac{1}{1-2\beta} \sum_{k=1}^n\bigg(\Big(\frac{1}{2}\Big)^{k}-\beta^{k}\bigg)= \frac{1}{1-\beta}.$$
As before, we deduce from Lemma~\ref{lemtoepbar} that
\begin{equation}
\label{cvgAn2}
\lim_{n \rightarrow + \infty}\frac{P_n}{|\dT_{n}|}=\frac{\sigma^2}{1-\beta}.
\hspace{1cm}\text{a.s.}
\end{equation}
Finally, Lemma~\ref{lemSumX} follows from the conjunction of (\ref{decosumX}), (\ref{resBn}), (\ref{resCn})
together with (\ref{cvgAn1}) and (\ref{cvgAn2}).\hspace{\stretch{1}} $\Box$\\

\noindent\textbf{Proof of Lemma~\ref{lemlimsumX} :} First of all, denote
\begin{equation*}
H_n = \sum_{k \in \dT_{n,p-1}}\dX_k
\hspace{1cm}\text{and}\hspace{1cm}
P_n = \sum_{k\in \dT_{n,p}}\veps_k,
\end{equation*}
As $|\dT_n|= 2^{n+1}-1$, we obtain from Equation~(\ref{defbarmatrix}) the recursive relation
\begin{eqnarray}
H_n & = & H_{p-1} + \sum_{k \in \dT_{n,p}}
\left(A_k\dX_{[\frac{k}{2}]} + \eta_k\right), \nonumber \\
& = & H_{p-1} + 2\overline{A}H_{n-1} +2\overline{a}(2^n-2^{p-1})e_1 +  P_ne_1
\label{limXone}
\end{eqnarray}
where $e_1=(1,0,\ldots,0)^t\in\dR^p$, $\overline{a}=(a_0+b_0)/2$ and the matrix
$$\overline{A}=\frac{A+B}{2}.$$
By induction, we deduce from (\ref{limXone}) that
\begin{eqnarray*}
\frac{H_n}{2^{n+1}} &\!=\!& \frac{H_{p-1}}{2^{n+1}} + \overline{A}\frac{H_{n-1}}{2^n} +\overline{a}\Bigl(\frac{2^{n}-2^{p-1}}{2^{n}}\Bigr)e_1 +  \frac{P_n}{2^{n+1}}e_1,  \\
&\!=\!&(\overline{A})^{n-p+1}\frac{H_{p-1}}{2^p} + 
\sum_{k=p}^n(\overline{A})^{n-k}
\left(\frac{H_{p-1}}{2^{k+1}}+ \overline{a}\Big(\frac{2^{k}-2^{p-1}}{2^{k}}\Big)e_1 + \frac{P_k}{2^{k+1}}e_1 \right).
\end{eqnarray*}
We have already seen via convergence (\ref{LLNeps1}) of Lemma~\ref{lemLLNeps} that
\begin{equation*}
\lim_{n \rightarrow + \infty}\frac{P_n}{2^{n+1}}=0
\hspace{1cm}\text{a.s.}
\end{equation*}
Finally, as $\|\overline{A}\|<1$, 
$$\sum_{n=0}^\infty\|(\overline{A})^n\|<\infty
\hspace{1cm} \text{and} \hspace{1cm} 
(\rI_{p}-\overline{A})^{-1}=\sum_{n=0}^\infty(\overline{A})^n,
$$ 
it follows from Lemma~\ref{lemtoepbar} that
\begin{equation*}
\lim_{n\rightarrow\infty}\frac{H_n}{2^{n+1}}=\overline{a}(\rI_{p}-\overline{A})^{-1}e_1
\hspace{1cm}\text{a.s.}
\end{equation*}
which ends the proof of Lemma~\ref{lemlimsumX}.\hspace{\stretch{1}}$ \Box$\\

\noindent\textbf{Proof of Lemma~\ref{lemlimsumXX} :} We shall proceed as in the proof of Lemma~\ref{lemlimsumX} and use the same notation. Let
\begin{equation*}
K_n = \sum_{k \in \dT_{n,p-1}}\dX_k\dX_k^t
\hspace{1cm}\text{and}\hspace{1cm} L_n = \sum_{k\in \dT_{n,p}}\veps_k^2.
\end{equation*}
We infer again from (\ref{defbarmatrix}) that
\begin{eqnarray*}
K_n&\!=\!& K_{p-1}+\sum_{k \in \dT_{n,p} }\left(A_k\dX_{[\frac{k}{2}] } +\eta_k\right)
\left(A_k\dX_{[\frac{k}{2}] } +\eta_k\right)^t\\
&\!=\!& K_{p-1}+ \sum_{k \in \dT_{n,p} }\!\! \veps^2_ke_1e_1^t+\sum_{k \in \dT_{n-1,p-1}}\!\!
\Bigl(A\dX_k\dX_k^tA^t + B \dX_k\dX_k^tB^t\Bigr)\\
\!\!\!\!&\!\!\!\!+\!& \sum_{k \in \dT_{n-1,p-1}}\!\!
\Bigl((a_0 +\veps_{2k})\dU_k(A)+(b_0 +\veps_{2k+1})\dU_k(B)
+2(\overline{a^2}+\zeta_{2k}) e_1e_1^t\Bigr)
\end{eqnarray*}
where $\dU_k(A)=A\dX_k e_1^t + e_1 \dX_k^t A^t$ and
$\dU_k(B)=B\dX_k e_1^t + e_1 \dX_k^t B^t$. In addition,
$\overline{a^2}=(a_0^2+b_0^2)/2$ and
$\zeta_{2k}=(a_0 \veps_{2k} + b_0\veps_{2k+1})$.
Therefore, we obtain that
\begin{equation}
\frac{K_n}{2^{n+1}} = \frac{1}{2} \left(A \frac{K_{n-1}}{2^n} A^t+B
\frac{K_{n-1}} {2^n}B^t\right)  + T_n  \label{limXXtwo}
\end{equation}
where
\begin{eqnarray*}
T_n &\!=\!& \left( \frac{L_n}{2^{n+1}}
+\overline{a^2}\Bigl( \frac{2^{n} -2^{p-1}}{2^{n}}\Bigr)
+\frac{1}{2^n}\sum_{k \in \dT_{n-1,p-1}}\!\!\zeta_{2k} \right)e_1e_1^t\\
&+& \frac{1}{2} \left( a_0\Bigl(A \frac{H_{n-1}}{2^n}e_1^t
+e_1\frac{H_{n-1}^t}{2^n}A^t\Bigr) + b_0\Bigl(B \frac{H_{n-1}}{2^n}e_1^t
+e_1\frac{H_{n-1}^t}{2^n}B^t\Bigr) \right) \\
&+& \frac{1}{2^{n+1}}\!\!
\sum_{k \in \dT_{n-1,p-1}}\!\!\Bigl(\veps_{2k}\dU_k(A) +
\veps_{2k+1}\dU_k(B)\Bigr).
\end{eqnarray*}
The two first results (\ref{LLNeps1}) and (\ref{LLNeps2}) of Lemma~\ref{lemLLNeps} together with Remark~\ref{rqCVeps2k} and Lemma~\ref{lem tree-generation} readily imply that
\begin{equation*}
\lim_{n \rightarrow + \infty}\frac{L_n}{2^{n+1}}=\sigma^2
\hspace{1cm}\text{a.s.}
\end{equation*}
and
\begin{equation*}
\lim_{n \rightarrow + \infty}
\frac{1}{2^n}
\sum_{k \in \dT_{n-1,p-1}}\!\!\zeta_{2k}=0 \hspace{1cm}\text{a.s.}
\end{equation*}
In addition, Lemma~\ref{lemlimsumX} gives
\begin{equation*}
\lim_{n \rightarrow + \infty}\frac{H_{n-1}}{2^{n}}=\lambda
\hspace{1cm}\text{a.s.}
\end{equation*}
Furthermore, denote 
\begin{equation*}
U_n=\sum_{k \in \dT_{n-1,p-1}}\!\!\Bigl(\veps_{2k}\dU_k(A)+
\veps_{2k+1}\dU_k(B)\Bigr).
\end{equation*}
For all $u\in\dR^p$, let $U_n(u)=u^tU_nu$. 
The sequence $\big(U_n(u)\big)$ is a real martingale transform.
Moreover, it follows from Lemma~\ref{lemSumX} that
$$
\sum_{k \in \dT_{n-1,p-1}}\!\!\Bigl|u^t\dU_k(A)u\Bigr|^2+\Bigl|u^t\dU_k(B)u\Bigr|^2
=\cO (|\dT_n|)
\hspace{1cm}\text{a.s.}
$$
Consequently, we deduce from the strong law of large numbers for martingale
transforms given in Theorem 1.3.24 of \cite{Duflo97} that
$U_n(u)=o(|\dT_n|)$ a.s. for all $u\in\dR^p$ which leads to $U_n=o(|\dT_n|)$ a.s. 
Therefore, we obtain that $(T_n)$ converges a.s. to $T$ given by
\begin{equation*}
T= ( \sigma^2 + \overline{a^2})e_1e_1^t+ \frac{1}{2} \left( A \lambda a_0 e_1^t
+ a_0 e_1 \lambda^tA^t + B \lambda b_0 e_1^t + b_0 e_1
\lambda^tB^t\right).
\end{equation*}
Finally, iteration of the recursive relation (\ref{limXXtwo}) yields 
\begin{equation*}
\frac{K_n}{2^{n+1}} = {\frac{1}{2^{n-p+1}}
\!\!\sum_{C \in \{A;B\}^{n-p+1}}\!\!C
\frac{K_{p-1}}{2^p}C^t}+\sum_{k=0}^{n-p} \,\frac{1}{2^k}\!\!
\sum_{C\in \{A;B\}^k}CT_{n-k}C^t.
\end{equation*}
On the one hand, the first term on the right-hand side converges a.s.
to zero as its norm is bounded $ \beta^{2(n-p+1)} \|K_{p-1}\|/2^p$.
On the other hand, thanks to Lemma~\ref{lempartX2}, the second term on 
the right-hand side converges to $\ell$ given by (\ref{eqrecl}), which completes
the proof of Lemma~\ref{lemlimsumXX}. .\hspace{\stretch{1}}$ \Box$\\

We now state a convergence result for the sum of $\|\dX_n\|^4$ which will be useful for the CLT.
\begin{Lemma}$\pt$
\label{lemsumX4}
Assume that $(\veps_n)$ satisfies \emph{(\textbf{H.1})} to \emph{(\textbf{H.5})}. Then, we have
\begin{equation}
\label{sumX4}
\sum_{k \in \dT_{n,p}}\|\dX_k\|^4 =\cO (|\dT_n|)
\hspace{1cm}\text{a.s.}
\end{equation}
\end{Lemma}

\noindent\textbf{Proof :} The proof is almost exactly the same as that of Lemma~\ref{lemSumX}. Instead of Equation~(\ref{decosumX}), we have
\begin{equation*}
\sum_{k \in \dT_{n,p}}  \|\dX_k\|^4  \leq \frac{64P_n}{(1-\beta)^3}+\frac{64\alpha^4Q_n}{(1-\beta)^3}+
8\overline{X}_1^4R_n
\end{equation*}
where
\begin{equation*}
P_n=\sum_{k \in \dT_{n,p}} \sum_{i=0}^{r_k-p}\beta^i\veps^4_{[\frac{k}{2^i}]},\quad
Q_n=\sum_{k \in \dT_{n,p}}\sum_{i=0}^{r_k-p}\beta^i,\quad
R_n=\sum_{k \in \dT_{n,p}}\!\!\beta^{4(r_k-p+1)}\,.
\end{equation*}
We already saw that $Q_n=\cO(|\dT_n|)$. In addition, it is not hard to see that $R_n=\cO(|\dT_n|)$. Therefore, we only need a sharper estimate for $u_n$. Via the same lines as in the proof of Lemma~\ref{lemSumX} together with the sharper results of Lemma~\ref{lemLLEps4TCL}, we can show that $P_n=\cO(|\dT_n|)$ a.s. which leads to (\ref{sumX4}).
\hspace{\stretch{1}}$ 
\Box$\\

\section{}
\begin{center}
{\bf On the quadratic strong law}
\end{center}
\label{appendixB}
We start with an auxiliary lemma closely related to the Riccation Equation for the inverse
of the matrix $S_n$.
\begin{Lemma}
\label{lemRiccati}
Let $h_n$ and $l_n$ be the two following symmetric square matrices of order $\delta_n$
$$h_n=\Phi_n^tS_n^{-1}\Phi_n
\hspace{1cm}\text{and}\hspace{1cm}
l_n=\Phi_n^tS_{n-1}^{-1}\Phi_n.
$$
Then, the inverse of $S_n$ may be recursively calculated as
\begin{equation}
\label{Riccaeq}
S_n^{-1} =S_{n-1}^{-1}-S_{n-1}^{-1}\Phi_{n}(\rI_{\delta_n}+l_n)^{-1}\Phi_n^tS_{n-1}^{-1}.
\end{equation}
In addition, we also have
$
(\rI_{\delta_n}-h_n)(\rI_{\delta_n}+l_n)=\rI_{\delta_n}.
$
\vspace{1ex}
\end{Lemma}

\begin{Remark}$\pt$
If $f_n=\Psi_n^t\Sigma_n^{-1}\Psi_n$, it follows from Lemma \ref{lemRiccati} that
\begin{equation}
\label{RiccaeqSig}
\Sigma_n^{-1} =\Sigma_{n-1}^{-1}-\Sigma_{n-1}^{-1}\Psi_{n}(\rI_{2\delta_n}-f_n)\Psi_n^t\Sigma_{n-1}^{-1}.
\end{equation}
\end{Remark}

\noindent\textbf{Proof :} As $S_n=S_{n-1}+\Phi_n\Phi_n^t$, relation~(\ref{Riccaeq})
immediately follows from Riccati Equation given e.g. in \cite{Duflo97} page 96.
By multiplying both side of (\ref{Riccaeq}) by $\Phi_n$, we obtain
\begin{eqnarray*}
S_n^{-1}\Phi_n &=& S_{n-1}^{-1}\Phi_n-S_{n-1}^{-1}\Phi_{n}(\rI_{\delta_n}+l_n)^{-1}l_n,\\
&=& S_{n-1}^{-1}\Phi_n-S_{n-1}^{-1}\Phi_{n}(\rI_{\delta_n}+l_n)^{-1}(\rI_{\delta_n}+l_n-\rI_{\delta_n}),\\
&=&S_{n-1}^{-1}\Phi_{n}(\rI_{\delta_n}+l_n)^{-1}.
\end{eqnarray*}
Consequently, multiplying this time on the left by $\Phi_n^t$, we obtain that
\begin{eqnarray*}
h_n &=& l_n(\rI_{\delta_n}+l_n)^{-1}=(l_n+\rI_{\delta_n}-\rI_{\delta_n})(\rI_{\delta_n}+l_n)^{-1},\\
&=&\rI_{\delta_n}-(\rI_{\delta_n}+l_n)^{-1}
\end{eqnarray*}
leading to $(\rI_{\delta_n}-h_n)(\rI_{\delta_n}+l_n)=\rI_{\delta_n}$. \hspace{\stretch{1}}$ \Box$\\

In order to establish the quadratic strong law for $(M_n)$, we are going to study separately the asymptotic behaviour of
$(\cW_n)$ and $(\cB_n)$ which appear in the main decomposition (\ref{maindecomart}).
\begin{Lemma}
\label{lemlimW}
Assume that $(\veps_n)$ satisfies \emph{(\textbf{H.1})} to \emph{(\textbf{H.3})}. Then, we have
\begin{equation}
\label{limW}
\lim_{n \rightarrow + \infty}\frac{1}{n} \cW_n= 2 \sigma^2
\hspace{1cm}\text{a.s.}
\end{equation}
\end{Lemma}

\noindent\textbf{Proof :} First of all, we have the decomposition
$\cW_{n+1}=\cT_{n+1}+\cR_{n+1}$ where
\begin{eqnarray*}
\cT_{n+1}&=&\sum_{k=p}^n \frac{\Delta M_{k+1}^t\Lambda^{-1}\Delta M_{k+1}}{|\dT_k|},\\
\cR_{n+1}&=&\sum_{k=p}^n \frac{\Delta M_{k+1}^t(|\dT_k|\Sigma_k^{-1}-\Lambda^{-1})\Delta M_{k+1}
}{|\dT_k|}.
\end{eqnarray*}
We claim that
\begin{equation*}
\lim_{n \rightarrow + \infty}\frac{1}{n} \cT_n= (p+1) \sigma^2
\hspace{1cm}\text{a.s.}
\end{equation*}
It will ensure via (\ref{cvgSigman}) that $\cR_n=o(n)$ a.s. leading to
(\ref{limW}). One can observe that $\cT_{n+1}=tr(\Lambda^{-1/2}H_{n+1}\Lambda^{-1/2})$ where
$$H_{n+1}=\sum_{k=p}^{n} \frac{\Delta M_{k+1}\Delta M_{k+1}^t}{|\dT_k|}.$$
Our goal is to make use of the strong law of large numbers for
martingale transforms, so we start by adding and
subtracting a term involving the conditional expectation of $\Delta
H_{n+1}$ given $\cF_n$. We have already seen in Section~\ref{sectionmartingale} that
for all $n\geq p-1$, $\dE[\Delta M_{n+1}\Delta M_{n+1}^t|\cF_n]=\Gamma
\otimes \Phi_n\Phi_n^t$. Consequently, we can split $H_{n+1}$ into
two terms
$$
H_{n+1}=\sum_{k=p}^{n} \frac{\Gamma \otimes \Phi_k\Phi_k^t}{|\dT_k|}+K_{n+1}
$$
where
$$
K_{n+1}=\sum_{k=p}^{n} \frac{\Delta M_{k+1}\Delta M_{k+1}^t-\Gamma \otimes \Phi_k\Phi_k^t}{|\dT_k|}.
$$
On the one hand, it follows from convergence (\ref{cvgSn})
and Lemma~\ref{lem tree-generation} that
\begin{equation*}
\lim_{n \rightarrow + \infty}\frac{\Phi_n\Phi_n^t}{|\dT_n|} = \frac{1}{2}L
\hspace{1cm}\text{a.s.}
\end{equation*}
Thus, Cesaro convergence yields
\begin{equation}
\label{cvgsumphi}
\lim_{n \rightarrow + \infty}\frac{1}{n}\sum_{k=p}^{n} \frac{\Gamma \otimes \Phi_k\Phi_k^t}{|\dT_k|}
= \frac{1}{2}(\Gamma \otimes L)
\hspace{1cm}\text{a.s.}
\end{equation}
On the other hand, the sequence $(K_n)$ is obviously a matrix martingale transform
satisfying
\begin{equation*}
\Delta K_{n+1}= K_{n+1}-K_n=\frac{1}{|\dT_{n+1}|}
\sum_{i,j\in\mathbb{G}_n}\Gamma_{ij}\otimes\left(
\begin{array}{ll}
1&\dX_j^t\\
\dX_i&\dX_i\dX_j^t
\end{array}
\right)
\end{equation*}
where
$$
\Gamma_{ij}=
\left(
\begin{array}{ll}
\varepsilon_{2i}\varepsilon_{2j}-\ind_{i=j}\sigma^2 \; & \varepsilon_{2i}\varepsilon_{2j+1}-\ind_{i=j}\rho\\
\varepsilon_{2i+1}\varepsilon_{2j}-\ind_{i=j}\rho \; &
\varepsilon_{2i+1}\varepsilon_{2j+1}-\ind_{i=j}\sigma^2
\end{array}
\right).
$$
For all $u\in \dR^{2(p+1)}$, let $K_n(u)=u^tK_n u$.
It follows from tedious but straightforward calculations, together with
(\ref{sumX2}), (\ref{sumX4}) and the strong law of large numbers for martingale
transforms given in Theorem 1.3.24 of \cite{Duflo97} that
$K_n(u)=o(n)$ a.s. for all $u\in\dR^{2(p+1)}$ leading to $K_n=o(n)$ a.s. 
Hence, we infer from (\ref{cvgsumphi}) that
\begin{equation}
\label{cvgHGL}
\lim_{n \rightarrow + \infty}\frac{1}{n}H_n= \frac{1}{2}(\Gamma \otimes L)
\hspace{1cm}\text{a.s.}
\end{equation}
Finally, we find from (\ref{cvgHGL}) that
\begin{eqnarray*}
\lim_{n \rightarrow + \infty}\frac{1}{n}\cT_n &=& \frac{1}{2}tr(\Lambda^{-1/2}(\Gamma \otimes L)\Lambda^{-1/2})
\hspace{1cm}\text{a.s.}\\
&=& \frac{1}{2}tr((\Gamma \otimes L)\Lambda^{-1})
\hspace{1cm}\text{a.s.}\\
&=& \frac{1}{2}tr(\Gamma \otimes \rI_{p+1})=(p+1)\sigma^2
\hspace{1cm}\text{a.s.}\\
\end{eqnarray*}
which completes the proof of Lemma \ref{lemlimW}
\hspace{\stretch{1}}$ \Box$\\

\begin{Lemma}
\label{lemlimB}
Assume that $(\veps_n)$ satisfies \emph{(\textbf{H.1})} to \emph{(\textbf{H.3})}. Then, we have 
\begin{equation*}
\cB_{n+1}=o(n)\hspace{1cm}\textrm{a.s.}
\end{equation*}
\end{Lemma}

\noindent\textbf{Proof :} Recall that
\begin{equation*}
\cB_{n+1}=2\sum_{k=p}^n M_k^t\Sigma_k^{-1}\Delta M_{k+1}=2\sum_{k=p}^n M_k^t\Sigma_k^{-1}\Psi_k\xi_{k+1}.
\end{equation*}
The sequence $(\cB_n)$ is a real martingale transform satisfying
\begin{equation*}
\Delta \cB_{n+1}=\cB_{n+1}- \cB_{n}=2M_n^t\Sigma_n^{-1}\Psi_n\xi_{n+1}.
\end{equation*}
Consequently, via the strong law of large numbers for martingale
transforms \cite{Duflo97}, we find that either $(\cB_n)$ converges a.s. or
$\cB_{n+1}=o(\nu_n)$ a.s. where
$$
\nu_n= \sum_{k=p}^n M_k^t\Sigma_k^{-1}\Psi_k\Psi_k^t\Sigma_k^{-1}M_k.
$$
However, for all $n\geq 2^{p-1}$, $\Psi_n\Psi_n^t=\rI_2\otimes\Phi_n\Phi_n^t$
which implies that
\begin{equation*}
\nu_n= \sum_{k=p}^n M_k^t\Sigma_k^{-1}(\rI_2\otimes\Phi_k\Phi_k^t)\Sigma_k^{-1}M_k
= \sum_{k=p}^n M_k^t(\rI_2\otimes S_k^{-1}\Phi_k\Phi_k^tS_k^{-1})M_k.
\end{equation*}
Furthermore, it follows from Lemma~\ref{lemRiccati} that
\begin{equation*}
S_{n-1}^{-1}-S_n^{-1}=S_n^{-1}\Phi_n(\rI_{\delta_n}+l_n)\Phi_n^tS_n^{-1}\geq S_n^{-1}\Phi_n\Phi_n^tS_n^{-1}
\end{equation*}
as the matrix $l_n$ is definite positive. Therefore, we obtain that
\begin{equation*}
\nu_n \leq \sum_{k=p}^nM_k^t(\Sigma_{k-1}^{-1}-\Sigma_k^{-1})M_k=\cA_n.
\end{equation*}
Finally, we deduce from the main decomposition (\ref{maindecomart}) that
\begin{equation*}
\cV_{n+1}+ \cA_n=o(\cA_n)+\cO(n)
\hspace{1cm}\text{a.s.}
\end{equation*}
leading to $\cV_{n+1}=\cO(n)$ and $\cA_{n}=\cO(n)$ a.s. as $\cV_{n+1}$ and $\cA_{n}$ are non-negative, which implies in turn that $\cB_n=o(n)$ a.s. completing
the proof of Lemma~\ref{lemlimB}.\hspace{\stretch{1}}$ \Box$\\

\noindent\textbf{Proof of Lemma~\ref{lem lim V+A} :} Convergence (\ref{cvgVA}) immediately
follows from (\ref{maindecomart}) together with Lemmas~\ref{lemlimW} and ~\ref{lemlimB}.
\hspace{\stretch{1}}$ \Box$\\

\section{}
\begin{center}
{\bf On Wei's Lemma}
\end{center}
\label{appendixC}
In order to prove (\ref{th Mn2}), we shall apply Wei's Lemma given in
\cite{Wei87} page 1672, to each entry of the vector-valued martingale
\begin{equation*}
M_n= \sum_{k=p}^n \sum_{i \in \dG_{k-1}} \left( \begin{array}{cccc}
\veps_{2i}  \\
\dX_i\veps_{2i} \\
\veps_{2i+1} \\
\dX_i\veps_{2i+1}
\end{array}\right).
\end{equation*}
We shall only carry out the proof for the first $(p+1)$ of $M_n$
inasmuch as the proof for the $(p+1)$ last components follows exactly the same lines.
Denote
$$
P_n=\sum_{k=p}^n \sum_{i \in \dG_{k-1}}\veps_{2i}
\hspace{1cm}\text{and}\hspace{1cm}
Q_n=\sum_{k=p}^n \sum_{i \in \dG_{k-1}}\dX_i\veps_{2i}.
$$
On the one hand, $P_n$ can be rewritten as 
${\displaystyle P_n=\sum_{k=p}^n\sqrt{|\dG_{k-1}|}v_k}$ 
where
$$
v_n=\frac{1}{\sqrt{|\dG_{n-1}|}}\sum_{i \in \dG_{n-1}}\veps_{2i}.
$$
We clearly have $\dE[v_{n+1}|\cF_n]=0$, $\dE[v_{n+1}^2|\cF_n]=\sigma^2$ a.s.
Moreover, it follows from (\textbf{H.1}) to (\textbf{H.3}) together with Cauchy-Schwarz inequality
that
\begin{eqnarray*}
\dE[v_{n+1}^4|\cF_n]&\!=\!&\frac{1}{|\dG_{n}|^2}\!\!\sum_{i \in \dG_{n}}\dE[\veps_{2i}^4|\cF_n]
+\frac{3}{|\dG_{n}|^2}\!\!\sum_{i \in \dG_{n}}\sum_{j\neq i}\dE[\veps_{2i}^2|\cF_n]\dE[\veps_{2j}^2|\cF_n]\\
&\leq & 3\sup_{i\in \dG_n}\dE[\veps_{2i}^4|\cF_n]
\hspace{0.5cm}\text{a.s.}
\end{eqnarray*}
which implies that $\sup \dE[v_{n+1}^4|\cF_n]< +\infty$ a.s.
Consequently, we deduce from Wei's Lemma that for all $\delta>1/2$,
\begin{equation*}
P_n^2=o(|\dT_{n-1}| n^\delta)
\hspace{1cm}\text{a.s.}
\end{equation*}
On the other hand, we also have 
${\displaystyle Q_n=\sum_{k=p}^n\sqrt{|\dG_{k-1}|}w_k}$ 
where
$$
w_n=\frac{1}{\sqrt{|\dG_{n-1}|}}\sum_{i \in \dG_{n-1}}\dX_i\veps_{2i}.
$$
It is not hard to see that $\dE[w_{n+1}|\cF_n]=0$ a.s.
Moreover, for all $1\leq k \leq p$, let $w_n(k)$ be the $k$th coordinate of the vector $w_n$.
It follows from (\textbf{H.1}) to (\textbf{H.3}) and 
Cauchy-Schwarz inequality that for all $1\leq k \leq p$,
\begin{eqnarray*}
\dE[w_{n+1}(k)^4|\cF_n]&\!\!\leq\!\!&\frac{1}{|\dG_{n}|^2}\!\!\sum_{i \in \dG_{n}}
\!\!X_{[\frac{i}{2^{k-1}}]}^4\dE[\veps_{2i}^4|\cF_n]
\!+\!\frac{3\sigma^4}{|\dG_{n}|^2}\!\!\sum_{i \in \dG_{n}}\sum_{j\neq i}\!X_{[\frac{i}{2^{k-1}}]}^2X_{[\frac{j}{2^{k-1}}]}^2\\
&\!\!\leq \!\!& 3\sup_{i\in \dG_n}\dE[\veps_{2i}^4|\cF_n] 
\left(\frac{1}{|\dG_{n}|}\!\!\sum_{i \in \dG_{n}}X_{[\frac{i}{2^{k-1}}]}^2\right)^2
\hspace{0.5cm}\text{a.s.}
\end{eqnarray*}
Hence, we obtain from Lemma~\ref{lemlimsumXX} that for all $1\leq k \leq p$, 
$\sup \dE[w_{n+1}(k)^4|\cF_n]< +\infty$ a.s.
Once again, we deduce from Wei's Lemma applied to each component of $Q_n$ that for all $\delta>1/2$,
\begin{equation*}
\|Q_n\|^2=o(|\dT_{n-1}| n^\delta)
\hspace{1cm}\text{a.s.}
\end{equation*}
which completes the proof of (\ref{th Mn2}).\hspace{\stretch{1}}$ \Box$\\

\section{}
\begin{center}
{\bf On the convergence of the covariance estimator}
\end{center}
\label{appendixD}
\noindent
It remains to prove that
\begin{equation*}
\lim_{n\rightarrow\infty}\frac{1}{n}
\!\sum_{k\in\dT_{n-1,p-1}}\!\!(\wh{\varepsilon}_{2k}
-\varepsilon_{2k})(\wh{\varepsilon}_{2k+1}-\varepsilon_{2k+1})
=\lim_{n\rightarrow\infty}\frac{R_n}{2n}= (p+1)\rho
\hspace{1cm}\text{a.s.}
\end{equation*}
where
$$
R_n=\sum_{k\in\mathbb{T}_{n-1,p-1}}(\wh{V}_k-{V}_k)^t\rJ_2(\wh{V}_k-{V}_k).
$$
It is not possible to make use of the previous convergence (\ref{cvgVchap})
because the matrix 
\begin{equation*}
\rJ_2=\left(\begin{array}{cc}0&1\\1&0\end{array}\right)
\end{equation*}
is not positive definite. Hence, it is necessary to rewrite our proofs. Denote
\begin{equation*}
\cV'_n=M_n^t\Sigma_{n-1}^{-1/2}(\rJ_2\otimes\rI_{p+1})\Sigma_{n-1}^{-1/2}M_n.
\end{equation*}As in the proof of Theorem~\ref{thmaptheta}, we have the decomposition
\begin{equation}
\label{decoprime}
\cV'_{n+1}+\cA'_n=\cV'_1+\cB'_{n+1}+\cW'_{n+1}
\end{equation}
where
\begin{eqnarray*}
\cA'_n&=&\sum_{k=p}^n M_k^{t}\big(\rJ_2\otimes(S_{k-1}^{-1}-S_k^{-1})\big)M_k,\\
\cB'_{n+1}&=&2\sum_{k=p}^n M_k^t(\rJ_2\otimes S_k^{-1})\Delta M_{k+1},\\
\cW'_{n+1}&=&\sum_{k=p}^n \Delta M_{k+1}^t(\rJ_2\otimes S_k^{-1})\Delta M_{k+1}.
\end{eqnarray*}
First of all, via the same lines as in Appendix B, we obtain that
\begin{eqnarray*}
\lim_{n \rightarrow + \infty}\frac{1}{n}\cW'_n &=& 
\frac{1}{2}tr((\rJ_2\otimes L^{-1})^{1/2}(\Gamma \otimes L)(\rJ_2\otimes L^{-1})^{1/2})
\hspace{1cm}\text{a.s.} \\
&=& \frac{1}{2}tr(\Gamma \rJ_2 \otimes \rI_{p+1})=(p+1)\rho
\hspace{1cm}\text{a.s.}
\end{eqnarray*}
Next, $(\cB'_n)$ is a real martingale transform satisfying 
$\cB'_{n+1}=o(n)$ a.s. Hence, we find the analogous of convergence (\ref{cvgVA})
\begin{equation}
\label{cvgVAprime}
\lim_{n \rightarrow + \infty}\frac{\cV'_{n+1}+\cA'_n}{n} = (p+1) \rho
\hspace{1cm}\text{a.s.}
\end{equation}
Furthermore, it follows from Wei's Lemma that for all $\delta>1/2$,
\begin{equation}
\label{cvgVprime}
\cV'_n=o(n^{\delta})\hspace{1cm}\text{a.s.}
\end{equation}
Therefore, we infer (\ref{decoprime}), (\ref{cvgVAprime}) and (\ref{cvgVprime}) that
\begin{equation}
\label{cvgAprime}
\lim_{n \rightarrow + \infty}\frac{1}{n}\cA'_n = (p+1)\rho
\hspace{1cm}\text{a.s.}
\end{equation}
Finally, by the same lines as in the proof of the first part of Theorem~\ref{thmapsigmarho}, we find that
\begin{equation*}
\lim_{n\rightarrow\infty}\frac{R_n}{n}=
2\lim_{n\rightarrow\infty}\frac{\cA'_n}{n}= 2(p+1)\rho\hspace{1cm}\text{a.s.}
\end{equation*}  
which completes the proof of convergence (\ref{cvgVJchap}).\hspace{\stretch{1}}$ \Box$\\

\noindent
{\bf Acknowledgements}
\ \vspace{2ex}\\
The authors would like to thanks the anonymous referees for their very careful reading of the
manuscript and for their suggestion to extend the first version of the paper
to asymmetric BAR($p$) processes.

\bibliographystyle{acm}
\bibliography{babar}

\begin{thebibliography}{10}

\bibitem{BaZh04}
{\sc Basawa, I.~V., and Zhou, J.}
\newblock Non-{G}aussian bifurcating models and quasi-likelihood estimation.
\newblock {\em J. Appl. Probab. 41A\/} (2004), 55--64.

\bibitem{CoSt86}
{\sc Cowan, R., and Staudte, R.~G.}
\newblock The bifurcating autoregressive model in cell lineage studies.
\newblock {\em Biometrics 42\/} (1986), 769--783.

\bibitem{DeMa08}
{\sc Delmas, J.-F., and Marsalle, L.}
\newblock Detection of cellular aging in a galton-watson process.
\newblock {\em arXiv, 0807.0749\/} (2008).

\bibitem{Duflo97}
{\sc Duflo, M.}
\newblock {\em Random iterative models}, vol.~34 of {\em Applications of
  Mathematics}.
\newblock Springer-Verlag, Berlin, 1997.

\bibitem{Guy07}
{\sc Guyon, J.}
\newblock Limit theorems for bifurcating {M}arkov chains. {A}pplication to the
  detection of cellular aging.
\newblock {\em Ann. Appl. Probab. 17}, 5-6 (2007), 1538--1569.

\bibitem{GBPSDT05}
{\sc Guyon, J., Bize, A., Paul, G., Stewart, E., Delmas, J.-F., and Tadd{\'e}i,
  F.}
\newblock Statistical study of cellular aging.
\newblock In {\em CEMRACS 2004---mathematics and applications to biology and
  medicine}, vol.~14 of {\em ESAIM Proc.} EDP Sci., Les Ulis, 2005,
  pp.~100--114 (electronic).

\bibitem{HaHe80}
{\sc Hall, P., and Heyde, C.~C.}
\newblock {\em Martingale limit theory and its application}.
\newblock Academic Press Inc., New York, 1980.
\newblock Probability and Mathematical Statistics.

\bibitem{Ham94}
{\sc Hamilton, J.~D.}
\newblock {\em Time series analysis}.
\newblock Princeton University Press, Princeton, NJ, 1994.

\bibitem{Hug96}
{\sc Huggins, R.~M.}
\newblock Robust inference for variance components models for single trees of
  cell lineage data.
\newblock {\em Ann. Statist. 24}, 3 (1996), 1145--1160.

\bibitem{HuBa99}
{\sc Huggins, R.~M., and Basawa, I.~V.}
\newblock Extensions of the bifurcating autoregressive model for cell lineage
  studies.
\newblock {\em J. Appl. Probab. 36}, 4 (1999), 1225--1233.

\bibitem{HuBa00}
{\sc Huggins, R.~M., and Basawa, I.~V.}
\newblock Inference for the extended bifurcating autoregressive model for cell
  lineage studies.
\newblock {\em Aust. N. Z. J. Stat. 42}, 4 (2000), 423--432.

\bibitem{HwBaYeo09}
{\sc Hwang, S.~Y., Basawa, I.~V., and Yeo, I.~K.}
\newblock Local asymptotic normality for bifurcating autoregressive processes
  and related asymptotic inference.
\newblock {\em Statistical Methodology 6\/} (2009), 61--69.

\bibitem{Wei87}
{\sc Wei, C.~Z.}
\newblock Adaptive prediction by least squares predictors in stochastic
  regression models with applications to time series.
\newblock {\em Ann. Statist. 15}, 4 (1987), 1667--1682.

\bibitem{ZhBa05a}
{\sc Zhou, J., and Basawa, I.~V.}
\newblock Least-squares estimation for bifurcating autoregressive processes.
\newblock {\em Statist. Probab. Lett. 74}, 1 (2005), 77--88.

\bibitem{ZhBa05b}
{\sc Zhou, J., and Basawa, I.~V.}
\newblock Maximum likelihood estimation for a first-order bifurcating
  autoregressive process with exponential errors.
\newblock {\em J. Time Ser. Anal. 26}, 6 (2005), 825--842.

\end{thebibliography}

\end{document}